\documentclass[a4paper,twoside,reqno,11pt]{elsarticle}

\usepackage[utf8]{inputenc}
\usepackage[english]{babel}

\usepackage{a4wide}
\usepackage{amsmath,amsfonts,amssymb,amsgen,amsbsy,amsthm,}
\usepackage{eucal,esint,dsfont,mathrsfs,MnSymbol}
\usepackage{stmaryrd} 
\usepackage{nomencl}
\usepackage{graphicx}
\usepackage{graphics}
\usepackage{float} 
\usepackage[section]{placeins} 
\usepackage{algorithm} 
\usepackage[noend]{algpseudocode} 
\usepackage{xcolor}
\usepackage{multirow}

\usepackage{subfig}
\usepackage{dirtytalk}
\usepackage{appendix}
\usepackage{color}
 \usepackage{ulem}

\addto\captionsenglish{}

\newcommand{\ud}{\mathrm{d}}

\newcommand{\ue}{\mathrm{e}}
\newcommand{\ui}{\mathrm{i}}

\newcommand{\half}{{\textstyle{\frac{1}{2}}}}

\newcommand{\be}{\begin{equation}}
\newcommand{\ee}{\end{equation}}
\newcommand{\eqdef}{\stackrel{\text{\tiny{def}}}{=}} 

\makeatletter
\newenvironment{breakablealgorithm}
{
	\begin{center}
		\refstepcounter{algorithm}
		\hrule height.8pt depth0pt \kern2pt
		\renewcommand{\caption}[2][\relax]{
			{\raggedright\textbf{\fname@algorithm~\thealgorithm} ##2\par}%
			\ifx\relax##1\relax 
			\addcontentsline{loa}{algorithm}{\protect\numberline{\thealgorithm}##2}%
			\else 
			\addcontentsline{loa}{algorithm}{\protect\numberline{\thealgorithm}##1}%
			\fi
			\kern2pt\hrule\kern2pt
		}
	}{
		\kern2pt\hrule\relax
	\end{center}
}
\makeatother

\begin{document}
	
	\begin{frontmatter}
		
		\title{Integrating factor techniques applied to the Schr{\"o}dinger-like
		equations. Comparison with Split-Step methods.}
		\author{Martino Lovisetto, Didier Clamond and Bruno Marcos}
		\address{Universit\'e C\^ote d'Azur, CNRS UMR 7351, LJAD, 
Parc Valrose, 06108 Nice cedex 2, France.}
		
\begin{abstract}
The nonlinear Schr{\"o}dinger and the Schr{\"o}dinger--Newton equations 
model many phenomena in various fields. Here, we 
perform an extensive numerical comparison between splitting methods (often
employed to numerically solve these equations) and the integrating factor technique, also called Lawson method. 
Indeed, the latter is known to perform very well for the nonlinear Schr{\"o}dinger
equation, but has not been thoroughly investigated for the Schr{\"o}dinger--Newton
equation. Comparisons are made in one and two spatial dimensions, exploring
different boundary conditions and parameters values. We show that for the short range
potential of the nonlinear Schr{\"o}dinger equation, the integrating factor technique
performs better than splitting algorithms, while, for the long range potential 
of the Schr{\"o}dinger--Newton equation, it depends on the particular
system considered.
\end{abstract}
\end{frontmatter}

\section{Introduction}
		
The nonlinear Schr\"odinger and the Schr\"odinger--Newton (also called 
Schr\"odinger--Poisson) equations describe a large number of phenomena 
in different physical domains. These equations are nonlinear variants of 
the Schr\"odinger one, which, in non-dimensional units, reads
\begin{equation} \label{def-schro}
\ui\,\partial_t\,\psi\ +\ \half\,\nabla^2 \psi\ -\ V\,\psi\ =\ 0,
\end{equation}
where $\psi$ is a function of space and time, $\nabla^2$ is the Laplace 
operator and $V$ is a function of $\psi$, space and time.

For the nonlinear Schr\"odinger equation (hereafter NLS), the {\it local\/} 
{\it nonlinear\/} potential is  
\begin{equation}
V\ =\ g\,|\psi|^2,
\end{equation}
where $g$ is a coupling constant. For $g>0$ the interaction is repulsive, 
while it is attractive for $g<0$. 
The NLS equation describes various physical phenomena, such as 
Bose--Einstein condensates \citep{dalfovo1999theory}, laser beams in some 
nonlinear optical media \citep{hasegawa1989optical}, water wave packets 
\citep{kharif2008rogue}, etc.

In the case of the Schr\"odinger--Newton (SN) equation, the potential is 
given by the Poisson equation 
\begin{equation}
\nabla^2\,V\ =\ g\, |\psi|^2,
\end{equation}
where $g$ is a coupling constant, the interaction being attractive if 
$g>0$ and repulsive if $g<0$. It is therefore {\it nonlinear\/} and {\it 
nonlocal}, giving rise to collective phenomena \citep{binney}, appearing 
for instance in optics  
\citep{dabby1968thermal,rotschild2005solitons,rotschild2006long}, 
Bose--Einstein condensates  
\citep{guzman2006gravitational}, astrophysics and cosmology 
\citep{hu2000fuzzy,paredes2016interference,marsh2019strong} 
and theories describing the quantum collapse of the wave function 
\citep{diosi2014gravitation,penrose1996gravity}. It is also used as a 
numerical model to perform  cosmological simulations in the semi-classical 
limit \citep{widrow1993using}. 

The SN equation takes a slightly different form when applied in cosmology 
\citep{peebles}. Here, due to the expansion of the universe, the 
Poisson equation is modified \citep{sikivie2015} as 
\begin{align} 
\nabla^2\,V\ =\ a^{-1}\,g\left(|\psi|^2-1\right), \label{def-poisson-cosmo}
\end{align}
where $a(t)$ is a scaling factor. The modification of the Poisson equation 
\eqref{def-poisson-cosmo} ensures that the potential is finite in an 
infinite universe. 
	
These NLS and SN equations are special cases of the 
Gross--Pitaevskii--Poisson (GPP) equation 
\begin{equation}
	\ui\,\partial_t\,\psi\ +\ \half\,\nabla^2\,\psi\ -\ V_1\,\psi\ -\ V_2\,\psi\ =\  0, \qquad 
	\nabla^2 V_1\ =\ g_1 \left|\,\psi\,\right|^2, \qquad 
	V_2\ =\ g_2 \left|\,\psi\,\right|^2.
\end{equation}
This equation appears in many fields, such as optics 
\citep{conti2003route,izdebskaya2016vortex}, 
Bose-Einstein condensates \citep{o2000bose} and cosmology 
(to simulate scalar field dark matter) 
\citep{suarez2017cosmological,chavanis2016collapse,bernal2006scalar}. 
	
In order to solve the above equations, except for very special cases, 
numerical methods must be used. Two families of 
temporal numerical schemes are commonly used to solve the Schr\"odinger 
equation with a nonlinear potential: the integrating factor technique 
(generally attached with a Runge--Kutta scheme) and the Split-Step method. 
In this paper, we present an extensive comparison between integrating factor 
methods and splitting algorithms, considering both accuracy and 
computational speed. Comparisons are made exploring different types 
of boundary conditions, in one and two spatial dimensions, with parameters 
ranging in values close to many regimes of physical interest.
The main reason for choosing these methods is that, in the literature, 
Split-Step solvers are commonly used to integrate both the SN and the 
NLS equations, while the integrating factor has been applied to integrate 
the NLS with very performing results 
\citep{bader2019efficient,besse2017high}. A natural question, which is also the main motivation of this work, is how the integrating factor technique performs when considering the long range interactions of the SN system instead of short range ones of the NLS.
We show that the integrating factor performs better than 
splitting algorithms for local interactions (such as the NLS). When a long 
range interaction (such as the one appearing in the SN equation) is considered, the 
relative performance between the integrating factor and splitting algorithms depends 
on the system.

The paper is organized as follows. In section \ref{sect-methods}, the methods of 
numerical time integration are described. Section \ref{sect-results} concerns detailed 
comparisons between Split-Step integrators (order 2, 4 and 6 with fixed time-step 
and order 4 with adaptive time-step) and standard algorithms with adaptive time-step 
belonging to the Runge--Kutta family \citep{dormand1980family,tsitouras2011runge,alexander1990solving} 
together with the integrating factor technique. 
Conclusions are drawn in section \ref{sect-conclusions}.

\section{Numerical algorithms}	\label{sect-methods}

In this section, we describe the different numerical methods used for 
the temporal resolution of the equations considered in the paper. First, 
we outline the integrating factor technique attached with an adaptive 
embedded scheme of arbitrary order. Then, we describe the Split-Step 
algorithm, with both fixed and adaptive time-step. 

Spatial resolution are performed with pseudo-spectral 
methods. In particular, we rely on methods based on 
fast Fourier transform ({FFT}) due to their efficiency 
and accuracy \citep{CanutoEtAl2006-1}.

\subsection{Integrating Factor}	\label{sect-IF}

The integrating factor method can be applied to any  differential equation of the form
\begin{equation}
\ui\,\partial_t\,\psi\ =\ F\!\left(\boldsymbol{r},t,\psi\right), \qquad
\psi\ =\ \psi\!\left(\boldsymbol{r},t\right),
\end{equation}
where the right-hand side can be split into linear and nonlinear parts. 
This results in
\begin{equation}    \label{eqIF1}
\ui\,\partial_t\,\psi\ +\ \mathcal{L}\,\psi\ =\ 
\mathcal{N}\!\left(\boldsymbol{r},t,\psi\right),
\end{equation}
where $\mathcal{L}$ is an easily computable (generally autonomous) linear  
operator and $\mathcal{N}\eqdef F+\mathcal{L}\psi$ is the 
remaining (usually) nonlinear part.
At the $n$-th time-step, with $t\in[t_{n},t_{n+1}]$, the 
change of dependent variable
\begin{equation}
\phi\ \eqdef\ \exp\!\left[\,(t-t_n)\,\mathcal{L}\,\right]
\psi \qquad\implies\qquad 
\ui\,\partial_t\,\psi\ =\ 
\exp\!\left[\,(t_n-t)\,\mathcal{L}\,\right]\left(\,\ui\,
\partial_t\,\phi\,-\,\mathcal{L}\,\phi\,\right),
\end{equation}
yields the equation  
\begin{equation} \label{eqIF2}
\ui\,\partial_t\,\phi\ =\ 
\exp\!\left[\,(t-t_n)\,\mathcal{L}\,\right]\mathcal{N}.
\end{equation}
Note that this change of variable is such that $\phi=\psi$ at $t=t_n$.
If the operator $\mathcal{L}$ is well chosen, the stiffness of 
\eqref{eqIF1} is considerably reduced and, for $t\in[t_n;t_{n+1}]$, 
the equation \eqref{eqIF2} can usually be well approximated by algebraic 
polynomials. Thus, standard time-stepping methods can efficiently 
solve \eqref{eqIF2}. Here, we focus on adaptive Runge--Kutta methods 
\cite{alexander1990solving,butcher2016numerical}.

As explained in \cite{lovisetto2022optimized}, it is possible to further improve and 
optimize the integrating factor method. The details of this improved version of the 
integrating factor technique are illustrated in \ref{AppIFC}. We perform our numerical 
tests using this optimized version of the integrating factor technique, which hereafter 
is denoted as IFC.

\subsubsection{Application to the Schr{\"o}dinger equation}

It is straightforward to apply these methods to the Schr{\"o}dinger 
equation \eqref{def-schro}. Specifically, for the SN and NLS equations, 
the Laplacian term is a linear operator while the potential $V$ is a 
nonlinear operator. Switching to Fourier space in position, the equation 
becomes
\begin{equation} \label{schro-fourier}
\ui\,\partial_t\, \widehat{\/\psi\/\/}\ -\ \half\,k^{2}\,\widehat{\/\psi\/\/}\ 
-\ \widehat{\,V\/\psi\,}\ =\ 0,
\end{equation}
where \say{hats} denote Fourier transforms of the 
underneath quantity and $k \eqdef |\boldsymbol{k}|$ is 
the wavenumber. 
Therefore, the system is now in a form where the application of the 
integrating factor technique is straightforward. With the change of 
variable  $\phi(\boldsymbol{k},t)= \widehat{\psi}(\boldsymbol{k},t)\,
\ue^{\ui{k^2}(t-t_0)/2}$, one obtains
\begin{equation} \label{SNIF1}
\partial_t\,\phi\ =\ -\/\ui\,\ue^{\ui{k^2}(t-t_0)/2}\, 
\widehat{\,V\/\psi\,}.
\end{equation}

In order to perform our numerical tests, we use the Dormand and Prince 
5(4) \cite{dormand1980family} and Tsitouras 5(4) \cite{tsitouras2011runge} 
integrators. 
Both schemes are Runge–Kutta pairs of order 5(4). However, we observe 
a speed difference between these solvers of maximum 10\%, depending 
on the simulated system. For this reason, we choose for each case the 
fastest of the two: specifically for NLS and the 
periodical SN equations we use the Dormand and Prince 
scheme, while in all the other cases we rely on Tsitouras' one.  
The higher-order Fehlberg 7(8) integrator \cite{alexander1990solving} is also 
used as {\it reference solutions} for accuracy comparisons (see section 
\ref{estimators}).

\subsection{Split-Step methods} \label{AdaptiveRK}

The Split-Step method \cite{blanes2008splitting} performs 
the temporal resolution of the Schr{\"o}dinger equation 
separating the linear terms from the nonlinear ones, in 
a different manner compared with the integrating factor. 
Writing the equation as 
\begin{equation}
\label{Schrodinger-H}
\ui\,\partial_t\,\psi\ =\ H\, \psi,
\end{equation}
$H=\half \nabla^2 + V$ being the Hamiltonian operator, the formal solution is
\begin{equation}
\label{formal-solution}
\psi(\boldsymbol{r}, t)\ =\ \exp\!\left({-\ui\int_{t_n}^t H\,\ud\/t}\right)\/\psi(\boldsymbol{r}, 
t_n), \qquad t\in[t_n;t_{n+1}].
\end{equation} 
Except for very few cases, the result of the operator $\exp({-\ui H(t-t_n)})$
applied to $\psi(\boldsymbol{r}, t_n)$ is unknown. 
Nevertheless, for $t\in[t_n;t_{n+1}]$, it is possible to approximate  
$\exp({-\ui H(t-t_n)})$ as a product of exponentials, each one involving either the potential or the Laplacian term, with appropriate coefficients. 
For example, the approximation corresponding to the Split-Step method 
of order 2 is
\begin{equation}\label{ss2tsf}
\ue^{-\int_{t_n}^t\ui\/H\,\ud\/t}\ =\ \ue^{-\ui K{(t-t_n)}/{2}}\,
\ue^{-\int_{t_n}^t\ui\/V\,\ud\/t}\,\ue^{-\ui 
K{(t-t_n)}/{2}}\ +\ O\left((t-t_n)^2\right),
\end{equation}
where $K=\half \nabla^2$.

At higher orders, the approximation of the operator $\exp({-\ui H(t-t_n)})$ 
is known as Suzuki-Trotter expansion \cite{suzuki1991general}. It is generally 
more complicated than \eqref{ss2tsf} and not unique, which can 
be determined with the Baker--Campbell--Hausdorff formula 
\cite{holden2010splitting}. 
For our numerical tests, we consider the Split-Step of orders 2, 4 and 
6, whose pseudo-codes are listed in \ref{app-ss}.

It is possible to design an adaptive time-step scheme with Split-Step 
methods. Here, we consider an adaptive embedded splitting pair 
\cite{thalhammer2012numerical} of order 4(3). This algorithm is 
characterized by a fourth-order splitting solver derived by Blanes and 
Moan \cite{blanes2002practical} embedded with third order scheme 
constructed by Thalhammer and Abhau 
\cite{thalhammer2012numerical,thalhammer2008high}.
The  pseudo-code for this algorithm, hereafter denoted ``SSa", is described 
in \ref{app-ss}.

\section{Numerical comparison of the different time-integrators}
\label{sect-results}

In this section, we compare the efficiencies of the methods previously 
described. The comparisons focus on speed and accuracy of each algorithm, 
simulating systems with different potentials, boundary conditions and 
physical regimes.
First, we outline the different estimators employed to determine the 
accuracy of each numerical integrator. Then, we list and summarize the 
results for every equations considered, in one and two spatial dimensions. 
We start with the NLS equation which is used as benchmark, since an 
analytical solution is known in the one dimensional case. Then we switch 
to the SN equation with both open and periodic boundary conditions. Finally, 
we present the results for the two dimensional Gross--Pitaevskii--Poisson 
equation, which can be considered as a hybrid version of the SN and NLS 
systems.

\subsection{Estimators of the accuracy of the time-integration algorithms}
\label{estimators}

The accuracy of each time-integration algorithm is estimated looking at 
three different indicators:
\begin{enumerate}
	\item The energy conservation. The energy $E$ is a constant of motion 
	for the Schr{\"o}dinger equation. For the NLS and SN equations, it is 
	defined as
	\begin{equation}
	E=\frac{1}{2}\int \ud\boldsymbol{r}\, \psi^*\,(-\nabla^2+V)\,\psi.
	\end{equation}
	Using the initial energy as reference, the error on the energy 
	conservation is
    \begin{equation}
    \label{deltaene}
    \Delta E_i = \left| \frac{E(t_i)}{E(t_0)} -1\right|,
    \end{equation} 
    where $t_0$ is the initial time and $t_i$ denotes the 
    $i$-th time-step of the numerical integration. This error being (in 
    general) time-dependent, we consider the error
    \begin{equation}
    \label{deltaene2}
    \overline{\Delta E} = \max_{i}\left[\Delta E_i\right],
    \end{equation}
    the latter being the maximal difference with respect to the initial value, during the whole simulation.
    
    \item Another constant of motion for the considered equation, is the mass, 
    \begin{equation}
	M=\int \ud\boldsymbol{r}\, |\psi|^2.
	\end{equation}
	This quantity is automatically conserved with machine precision when using splitting algorithms, while it is not in general the case with the integrating factor. For this reason, when the latter technique is employed, we impose mass conservation at each time-step, multiplying the solution $\psi$ by $M_0/\int \ud\boldsymbol{r}\, |\psi|^2$, where $M_0$ is the initial mass.

	\item The error on the solution performing time reversion tests. This 
	quantity is obtained running a simulation up to a given time 
	$t_{\mathrm{fin}}$, then reversing the time and evolving back to the 
	initial instant. The error is monitored using the $L_{\infty}$-norm of 
	the difference between the solution at the initial time, at beginning of 
	the simulation and at the end of it. Denoting the \say{backward} solution 
	by $\Delta \psi_{\mathrm{rev}}$, one has
	\begin{equation}
	\label{deltapsi}
	\Delta \psi_{\mathrm{rev}} = \max_{i}\left(\,\left|\, \psi(x_i,
	t_0)\, -\, \psi_{\mathrm{backward}}(x_i,t_0)\,\right|\,\right).
	\end{equation}
	
	\item The two estimators above favorize {\it a priori} time-splitting 
	algorithms because they are symplectic and reversible, whereas the 
	integrating factor is not. For this reason, we also compare the result 
	of the simulations with a ``reference one'', very accurate, using an 
	adaptive Fehlberg integrator of order $7$ embedded within an order $8$ 
	scheme, with a very small tolerance, $\mathrm{tol}=10^{-14}$. Defining 
	this estimator as $\Delta \psi_{\mathrm{ref}}$, one has
	\begin{equation}
	\label{deltapsiref}
	\Delta \psi_{\mathrm{ref}} = \max_{i}\left(\,\left|\,\psi(x_i,t_{f})\, 
	-\, \psi_{\mathrm{F7(8)}}(x_i,t_{\mathrm{fin}})\,\right|\,\right),
	\end{equation}
	where $\psi$ is the numerical solution provided by the particular 
	method considered and $\psi_{\mathrm{F7(8)}}$ is the one outputted by the 
	Fehlberg 7(8) integrator.
\end{enumerate}

\subsection{1D nonlinear Schr{\"o}dinger equation}

We first consider the case of the one dimensional NLS 
\begin{equation}
\label{nls}
\ui \frac{\partial \psi}{\partial t}\ +\ \frac{1}{2}\, \frac{\partial^2 
\psi}{\partial x^2}\ +\, \left|\psi \right|^2\, \psi\ =\ 0,
\end{equation}
which admits a  simple analytical solution
\begin{equation}
\psi(x,t)\ =\ \sqrt{2}\,{\mathrm{sech}}\left(\sqrt{2}\,x\right)
\/\exp({\ui \,t}).
\end{equation}
We present a set of simulations in order to compare the Split-Step 
integrators with the IFC, looking at the energy conservation error, the 
error on the solution and the total time needed to run each simulation. 
In these simulations, the space is discretized with $N=2048$ points, in a 
domain of length $L=80$ and the analytical solution at $t=0$ is used as 
initial condition.
The results are summarized in table (\ref{nlst3}). 
We observe that the IFC solver is the fastest one by at least a factor 2, 
presenting at the same time the best results to all the indicators: it uses 
a larger time-step, presents equal or better energy conservation, returns 
only a slightly worse $\Delta \psi_{\mathrm{rev}}$ and it is one of the 
best comparing to the reference simulation. 

\begin{table}
	\centering
	\begin{tabular}{ |c|c|c|c|c|c| } 
		\hline
		Method & $\Delta t$ & $\overline{\Delta E}$ & $\Delta \psi_{\mathrm{rev}}$ & $\Delta \psi_{\mathrm{ref}}$ & $T(s)$ \\
		\hline
		SS2 & $ 10^{-3}$& $7.5 \cdot 10^{-13}$ & $6.3 \cdot 10^{-10}$ & $1.1 \cdot 10^{-4}$ & 86.7\\ 
		\hline
		SS4 & $5 \cdot 10^{-3}$& $ 10^{-15}$& $8.2 \cdot 10^{-10}$ & $5.8 \cdot 10^{-7}$ & 38.1\\ 
		\hline
		SS6 & 2 $\cdot 10^{-2}$& $ 10^{-15}$ & $2.0 \cdot 10^{-10}$ & $3.7 \cdot 10^{-9}$ & 29.3\\ 
		\hline
		SSa, $\mathrm{tol}=10^{-6}$ & $2.1 \cdot 10^{-2}$ & $1.6 \cdot 10^{-12}$  & $3.5 \cdot 10^{-10}$ & $6.5 \cdot 10^{-10}$ & 34.0\\ 
		\hline
		IFC, $\mathrm{tol}=10^{-9}$ & $2.1 \cdot 10^{-2}$ & $ 10^{-15}$ & $1.2 \cdot 10^{-9}$ & $7.8 \cdot 10^{-10}$ & 14.5 \\ 
		\hline
	\end{tabular}
	\caption{Comparison for the 1D NLS equation between the IFC method and the Split-Step solvers. $T$ is the total time required to run each simulation, measured in seconds. The $\Delta t$ for adaptive algorithms is the averaged one.}
	\label{nlst3}
\end{table}

\subsection{1D Schr{\"o}dinger--Newton equation \label{SN1Dsect}} 
\label{sect-1D-SN}
\begin{figure}
\includegraphics[width=0.32\textwidth]{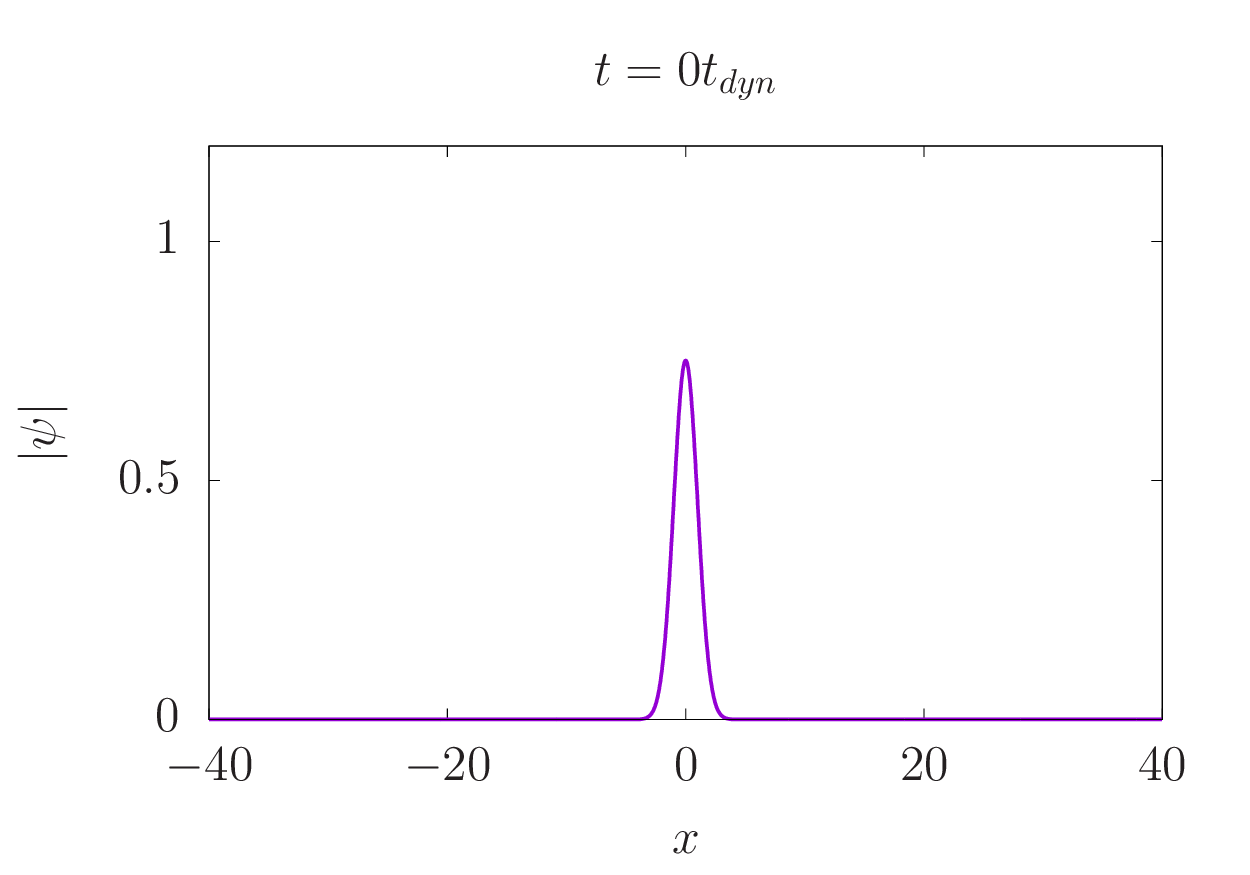}
\includegraphics[width=0.32\textwidth]{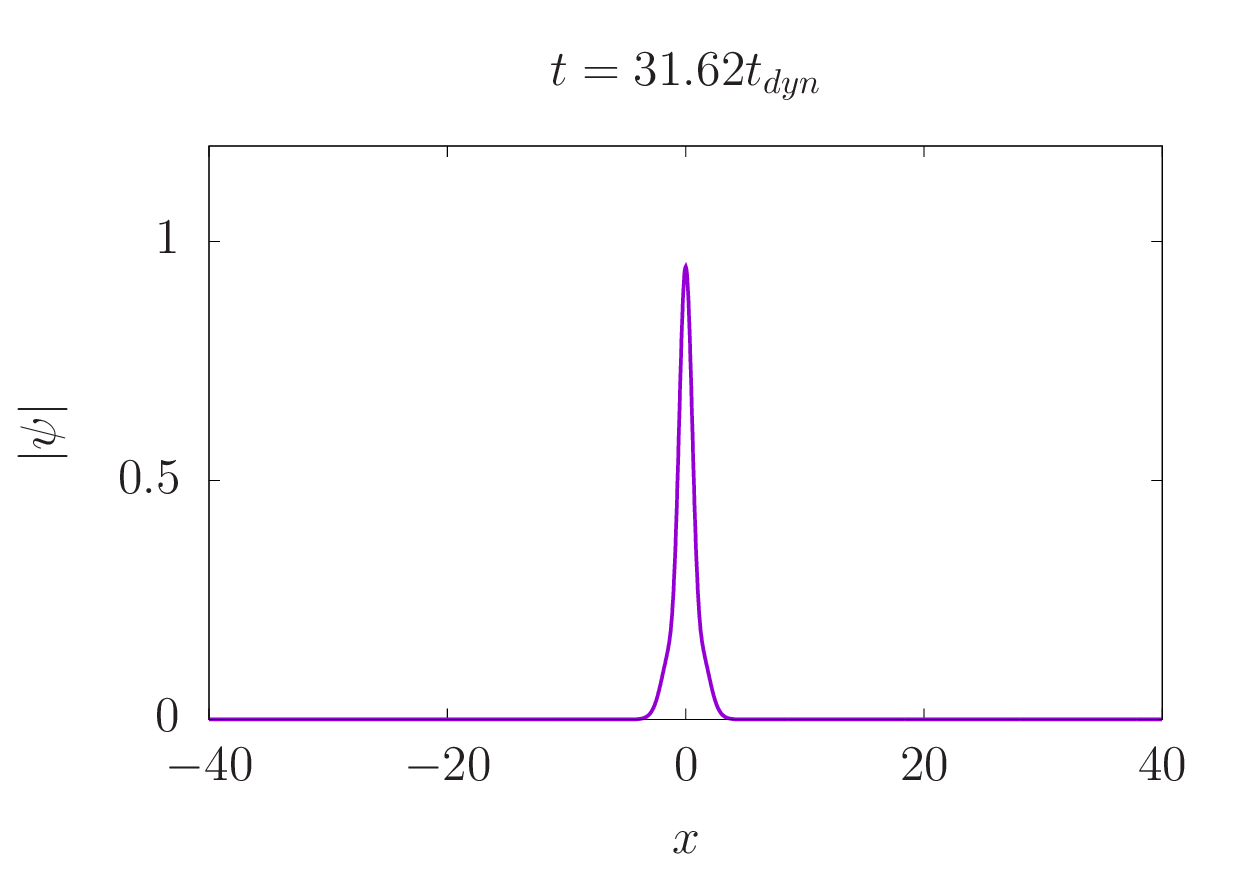}
\includegraphics[width=0.32\textwidth]{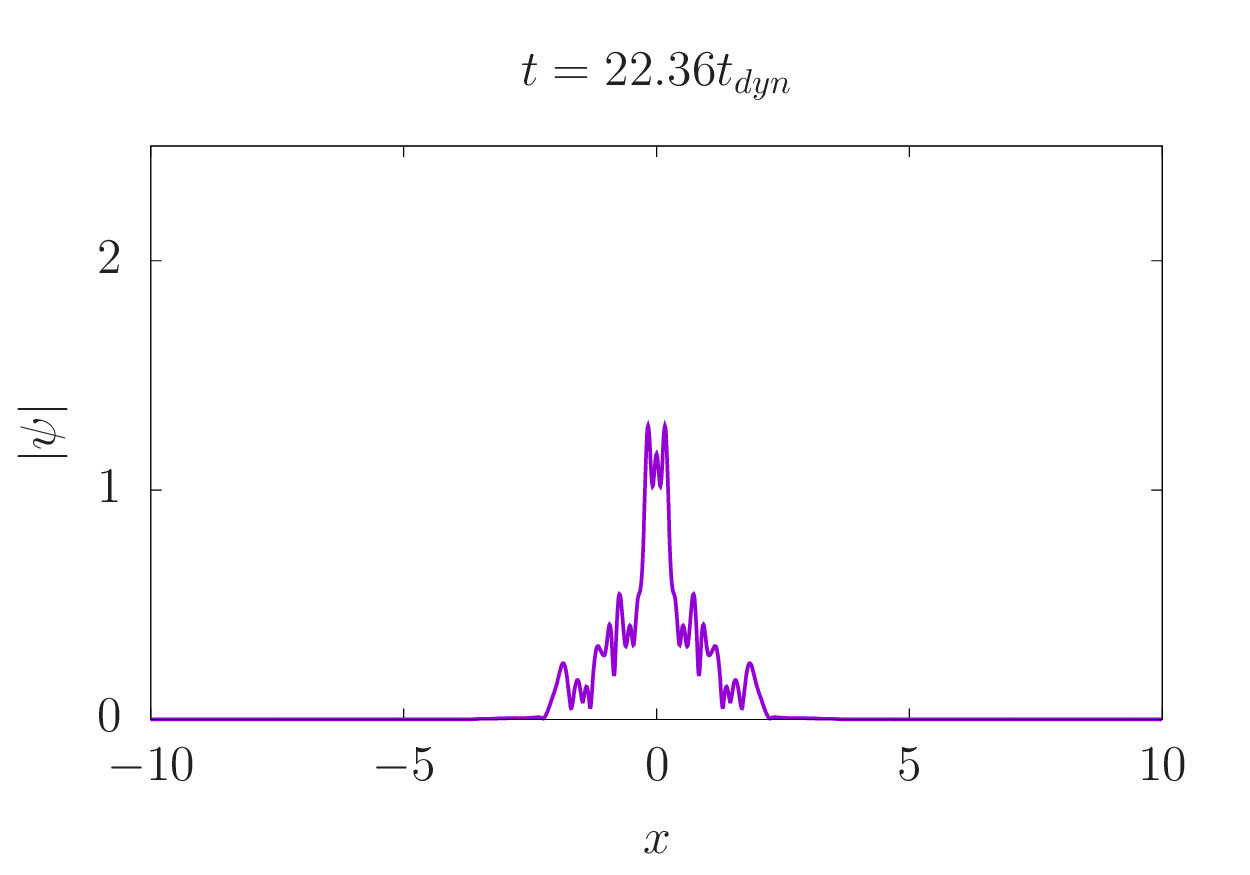}
	
	\caption{Snapshots of the modulus of the solution of the 1D SN equation $\left| \psi \right|$. The left plot is the initial condition, center plot and right plot correspond to the solution at the end of the simulation for $g=10$ and $g=500$ respectively.}
	\label{solsn1d}
\end{figure}

We now focus on the SN system, starting from the case of a single spatial 
dimension, 
\begin{equation}
\begin{gathered}
\label{sp}
\ui\, \frac{\partial\,\psi}{\partial t}\ +\ \frac{1}{2}\,\frac{\partial^2
\,\psi}{\partial x^2}\ -\  V\, \psi\ =\ 0, \qquad
\frac{\partial^2\,V}{\partial x^2}\ =\ g \left| \psi \right |^{2}.
\end{gathered}
\end{equation}
The solutions of \eqref{sp} depend on the initial condition and on the single 
parameter $g$. The chosen initial condition is $\psi(x,t=0)=
\exp({-x^2/2})/\sqrt[4]{\pi}$. The potential $V$ is calculated using 
Hockney's method \cite{hockney2021computer}. We perform a set of tests with 
different values of the parameter $g$, corresponding to different physical 
regimes. The case $g=10$ corresponds to a system in the quantum regime, 
i.e., with an associated De Broglie wavelength of the order of the size of 
the system, while $g=500$ corresponds to a system closer to the semi-classical 
regime, with an associated De Broglie wavelength about 20 times smaller 
than the size of the system. The typical evolution of this system is characterized by the initial condition which oscillates, exhibiting a 
complex dynamics. This is particularly visible in the semi-classical regime, 
in which high frequency oscillations appear in the wavefunction, as shown in Fig.~\ref{solsn1d}.
The simulation is run in a domain of length $L=80$, discretized into 
$N=2048$ points in the $g=10$ case, while for $g=500$ we set $L=20$ and 
$N=2048$. The characteristic time of dynamics is defined as 
$t_{\mathrm{dyn}}=\left|g\right|^{-1/2}$.
\begin{figure}
	\minipage{0.47\textwidth}
	\resizebox{\textwidth}{!}{\includegraphics{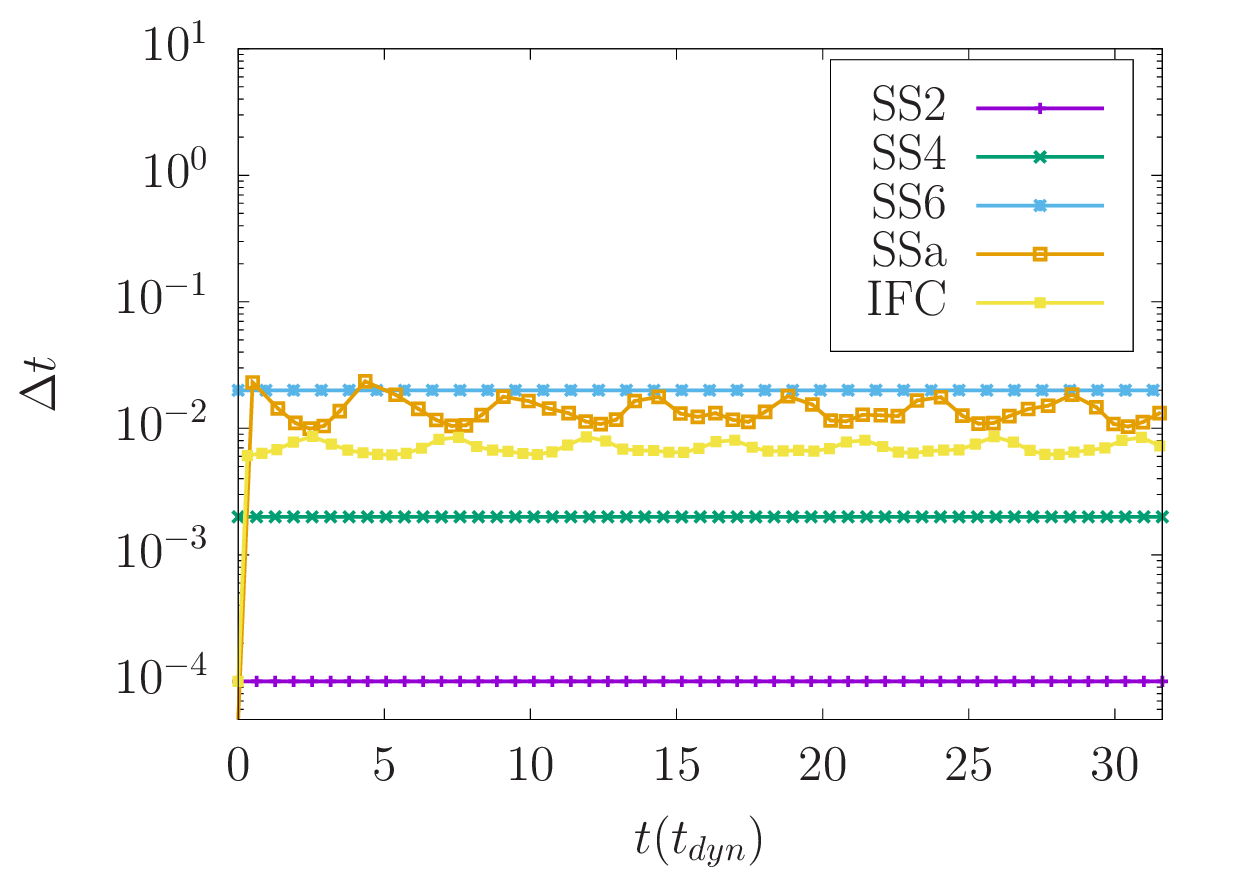}}
	\endminipage\hfill
	\minipage{0.47\textwidth}%
	\resizebox{\textwidth}{!}{\includegraphics{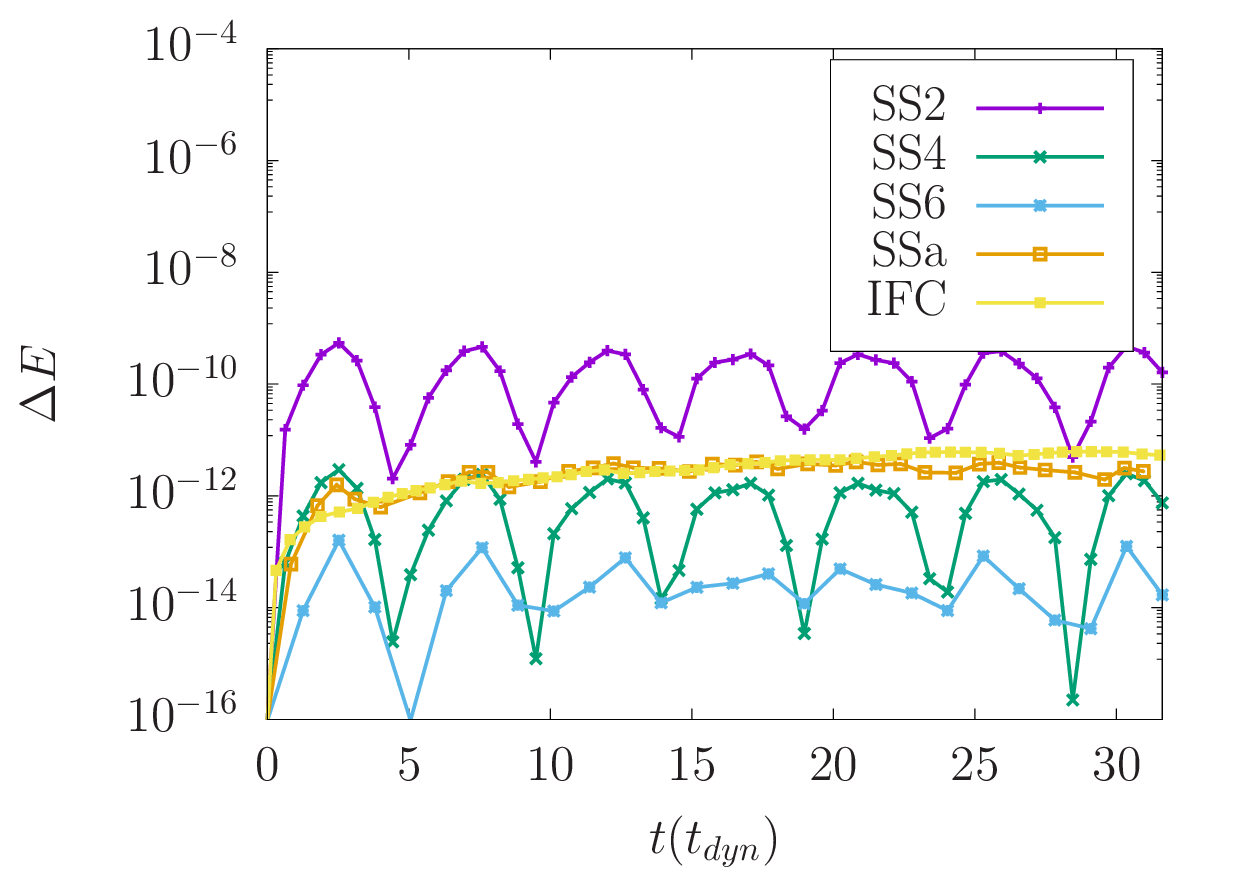}} 
	\endminipage
	
	\minipage{0.47\textwidth}
	\resizebox{\textwidth}{!}{\includegraphics{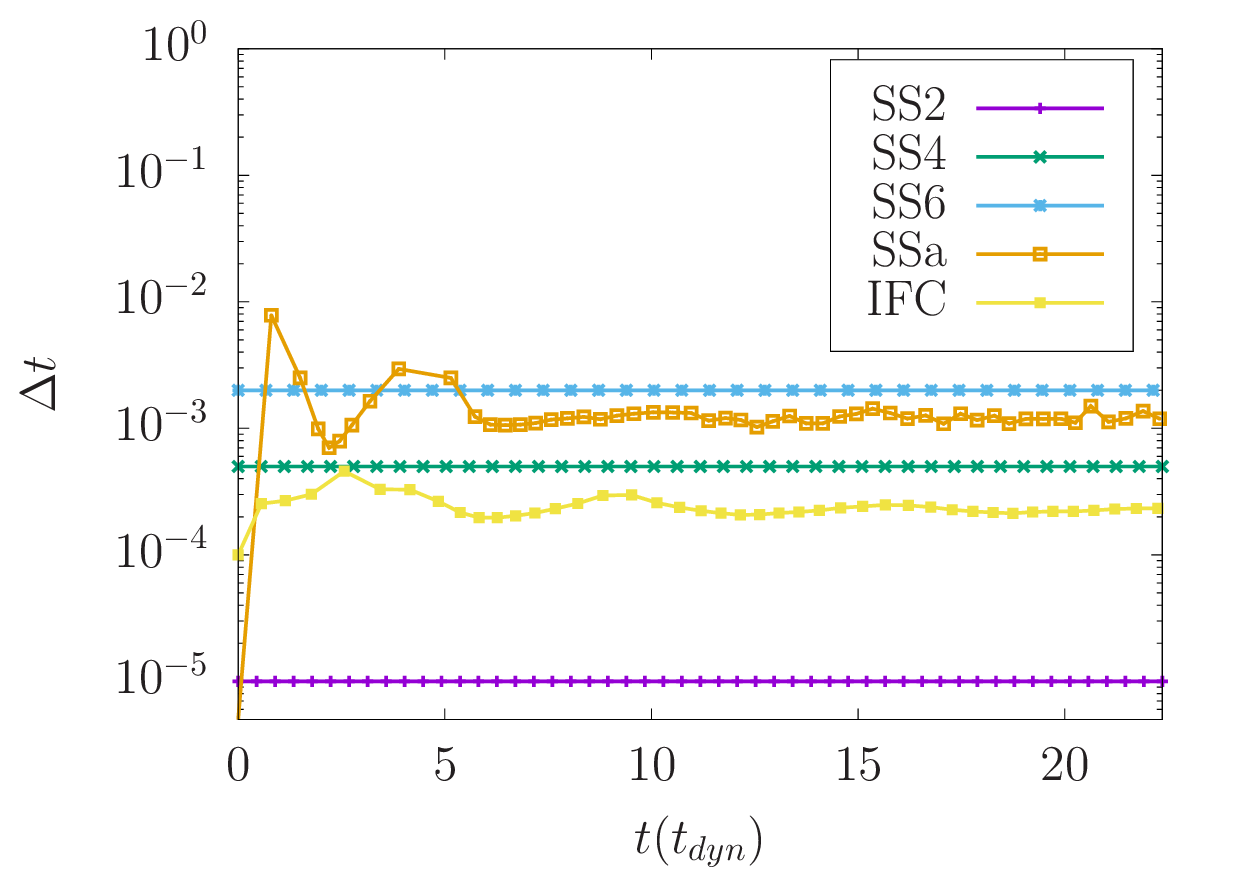}}
	\endminipage\hfill
	\minipage{0.47\textwidth}%
	\resizebox{\textwidth}{!}{\includegraphics{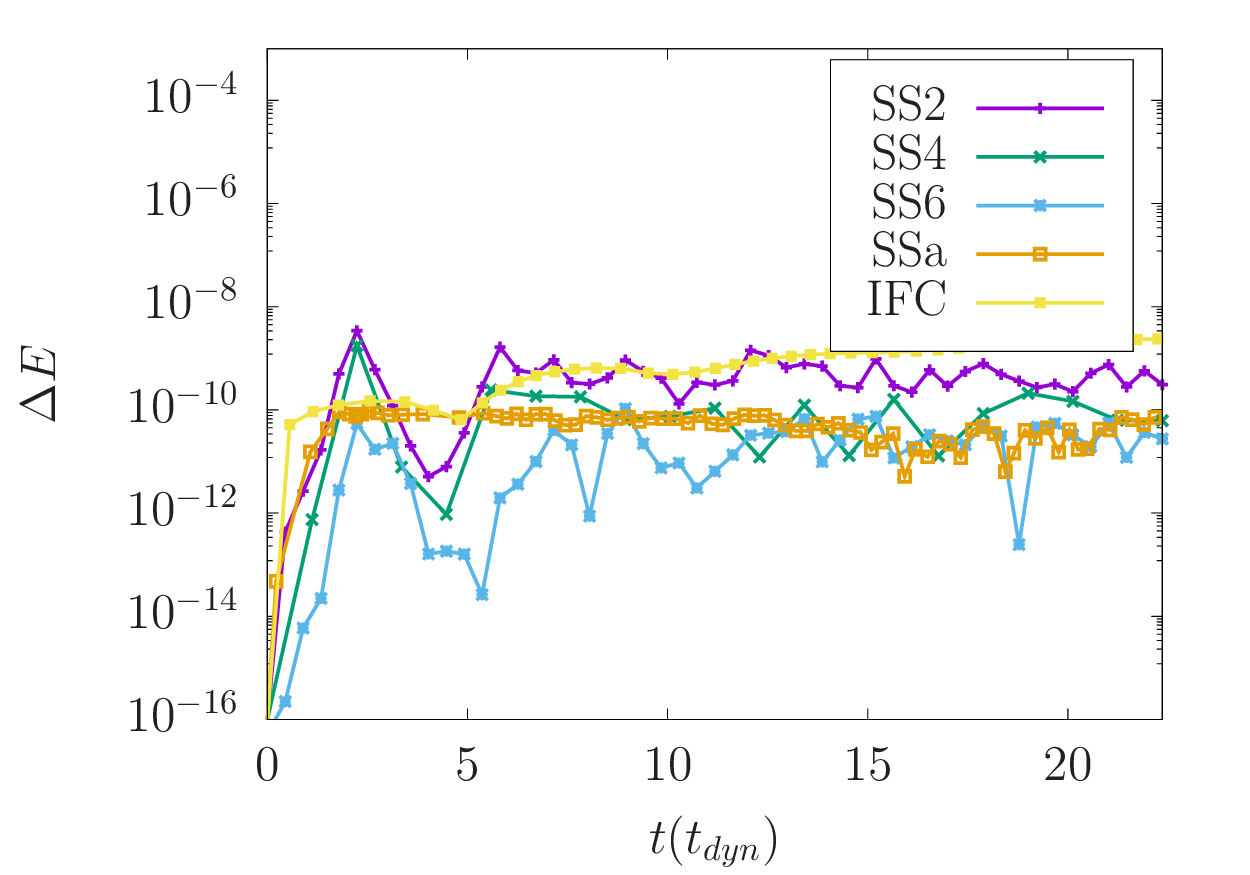}} 
	\endminipage
	\caption{Comparison for the 1D SN equation between the time-step and the error on the energy conservation  for the IFC method and the Split-Step solvers for both the cases $g=10$ (upper plots) and $g=500$ (lower plots). }
	\label{sp1df1b}
\end{figure}
\begin{table}
	\centering
	\begin{tabular}{ |c|c|c|c|c|c|c| } 
		\hline
		g& Method & $\Delta t/t_{\mathrm{dyn}}$ & $\overline{\Delta E}$ & $\Delta \psi_{\mathrm{rev}}$ & $\Delta \psi_{\mathrm{ref}}$ & $T(s)$ \\
		\hline
		\multirow{4}{3em}{10}&SS2 & $ 3.1 \cdot 10^{-4}$& $5.5 \cdot 10^{-10}$ & $4.5 \cdot 10^{-11}$ & $4.5 \cdot 10^{-9}$ & 123.3\\ 
		&SS4 & $6.3 \cdot 10^{-3}$&  $ 2.9 \cdot 10^{-12}$& $2.5 \cdot 10^{-12}$ & $5.5 \cdot 10^{-11}$ & 11.1\\ 
		&SS6 & $6.3 \cdot 10^{-2}$&  $ 1.5 \cdot 10^{-13}$& $4.0 \cdot 10^{-12}$ & $1.1 \cdot 10^{-11}$ & 3.5\\ 
		&SSa & $ 4.3\cdot 10^{-2}$& $ 4.2\cdot 10^{-12}$  & $ 8.5\cdot 10^{-11}$ & $ 5.9\cdot 10^{-11}$ & 6.4 \\ 
		&IFC &2.2 $\cdot 10^{-2}$& $ 1.8 \cdot 10^{-12}$ & $1.3 \cdot 10^{-10}$ & $3.3 \cdot 10^{-11}$ & 5.0 \\ 
		\hline
		\multirow{4}{3em}{500}&SS2 & $ 10^{-5}$& $6.0 \cdot 10^{-9}$ & $1.1 \cdot 10^{-10}$ & $5.2 \cdot 10^{-8}$ &120.0\\ 
		&SS4 & $5 \cdot 10^{-4}$&  $ 2.5 \cdot 10^{-9}$& $4.9 \cdot 10^{-12}$ & $1.3 \cdot 10^{-7}$ &5.4\\ 
		&SS6 & $2 \cdot 10^{-3}$&  $ 1.8 \cdot 10^{-10}$& $2.8 \cdot 10^{-11}$ & $7.6 \cdot 10^{-8}$ & 3.5\\
		&SSa & $ 3.2 \cdot 10^{-2}$& $1.0 \cdot 10^{-10}$ &$6.4 \cdot 10^{-8}$  & $2.9 \cdot 10^{-8}$ & 6.7\\ 
		&IFC&$5.4 \cdot 10^{-3}$& $1.8 \cdot 10^{-9}$ & $2.3 \cdot 10^{-8}$ & $1.5 \cdot 10^{-8}$ &11.9\\ 
		\hline
	\end{tabular}
	\caption{Comparison for the 1D SN equation between the IFC method and the Split-Step solvers.  The SSa simulations and IFC have been performed with a tolerance $\mathrm{tol}=10^{-7}$ and $\mathrm{tol}=10^{-10}$ respectively for $g=10$ and $\mathrm{tol}=10^{-6}$ and $\mathrm{tol}=10^{-10}$ respectively for $g=500$. The $\Delta t$ for adaptive algorithms is the averaged one. $T$ is the total time required to run each simulation, measured in seconds.}
	\label{sp1dt2}
\end{table}
In table (\ref{sp1dt2}), we compare the Split-Step integrators with the 
IFC, looking at the energy conservation error, the error on the solution and 
the total time needed to run each simulation. Here, splitting methods proved 
to be faster than the integrating factor. In addition, the SS4
and SS6 performed better than the adaptive integrators. This is 
due the fact that, for this particular system, the extra computational cost 
due to the implementation of the adaptive-step is not fully compensated by 
the time-gain in terms of computational speed.
Indeed, splitting algorithms with fixed time-step require a smaller number 
of computational operations to be implemented. For this reason, here, 
choosing a \say{proper} fixed time-step still results in a slightly faster 
numerical integration compared to an adaptive scheme.

\subsubsection{Periodical case}

We now switch to another version of the SN system, which has important 
applications in cosmology in order to simulate the formation of large-scale 
structures in the universe \eqref{def-poisson-cosmo}. We take $a=1$, which in cosmology corresponds to the case of a static universe \citep{weinberg2008cosmology}; we do not expect 
modifications of our conclusions for different cosmological models. In one 
dimension the equations read
\begin{subequations}
	\label{sp1dcosmo}
	\begin{align}
	\label{sp1dcosmo1}
	&\ui\, \frac{\partial\, \psi}{\partial t}\ +\ \frac{1}{2}\, 
	\frac{\partial^2\, \psi}{\partial x^2}\ -\  V\, \psi\ =\ 0, \\
	\label{sp1dcosmo2}
	&   \frac{\partial^2\,V}{\partial\/x^2}\ =\ g\,( \left| \psi \right |^{2}-1),
	\end{align}
\end{subequations}
where the wavefunction $ \left| \psi \right |^{2}$ is normalized to unity. 
The potential $V$ is obtained calculating the inverse of the Laplacian in 
Fourier space and transforming back the result to real space. We take 
``{cold}'' initial conditions (see \cite{davies1997,mocz2018}), namely,
\begin{equation}
\label{sp1dcosmoIC1}
\psi(x,t=0)\ =\ \sqrt{\rho_0+\delta \rho(x)}\,\exp({\ui \theta(x)}),
\end{equation}
where $\theta$ is a function whose gradient is proportional to the initial 
velocity field (set to zero for simplicity), $\rho_0$ is the 
background constant density and $\delta \rho (x)$ is the density fluctuations, 
generated as
\begin{equation}
\delta \rho(x)\ =\ \mathcal{F}^{-1}\left[\,R(k)\,\sqrt{P(k)}\,\right],
\end{equation}
where $R(k)$ is a Gaussian random field, with zero average and unity 
variance. The function $P(k)$ is called {\it Power Spectrum}, and it is 
defined as
\begin{equation}
P(k)\ =\ \frac{1}{L^d}\,\left|\,\widehat{\delta\rho}(k)\,\right|^2,
\end{equation}
corresponding to the initial density fluctuations one wants to generate. 
The initial conditions are numerically initialized applying an additional 
filter $F(k)$ in Fourier space with the aim of setting to zero all the 
modes corresponding to a space scale comparable (or smaller) than  
the grid-step
\begin{equation}
F(k)\ =\ \mathrm{sech}\left(\left(k/k_F\right)^{10}\right),
\end{equation}
with $k_F=k_{N}/8$, where $k_{N}$ is the Nyquist wavelength, defined 
as $k_{N}=\frac{N}{2L}$. Thus, the initial condition is
\be	
\psi(x,t=0)\ =\ \mathcal{F}^{-1}\left[\,F(k)\,\mathcal{F}\left[
\sqrt{\rho_0+\delta \rho(x)}\right]\,\right].
\ee

In the simulations, space is discretized with $N=1024$ points, 
in a domain of length $L=1$ and a constant power spectrum is used as 
initial condition. We show the simulation results in the semi-classical 
regime. The latter corresponds to large values of the parameter $g$, as 
one has $g \propto \hbar^{-2}$. Specifically, we take $g=10^6$ (we do 
not observe differences in the performance in the quantum regime, i.e., 
for smaller values of $g$) and $\rho_0=1$. 

\begin{figure}
\includegraphics[width=0.32\textwidth]{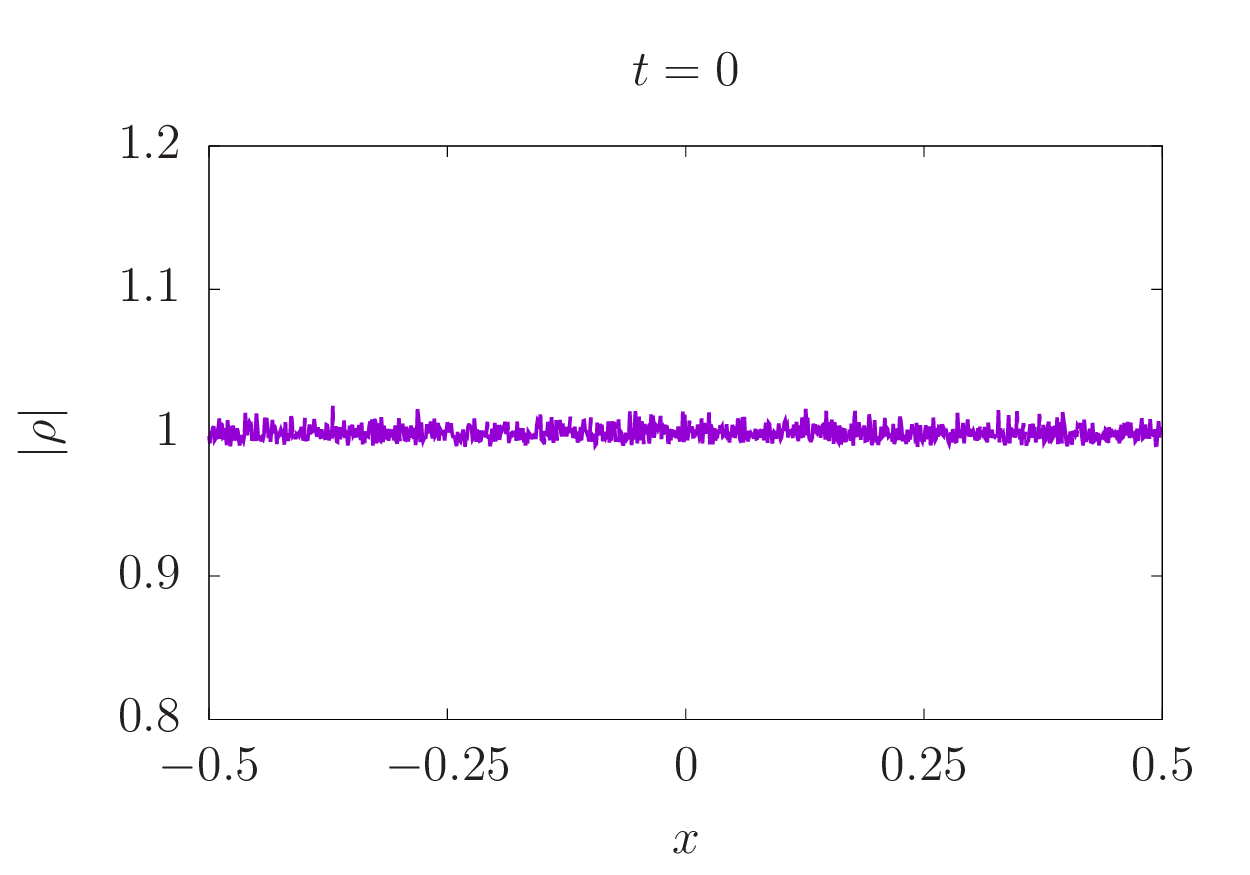}
\includegraphics[width=0.32\textwidth]{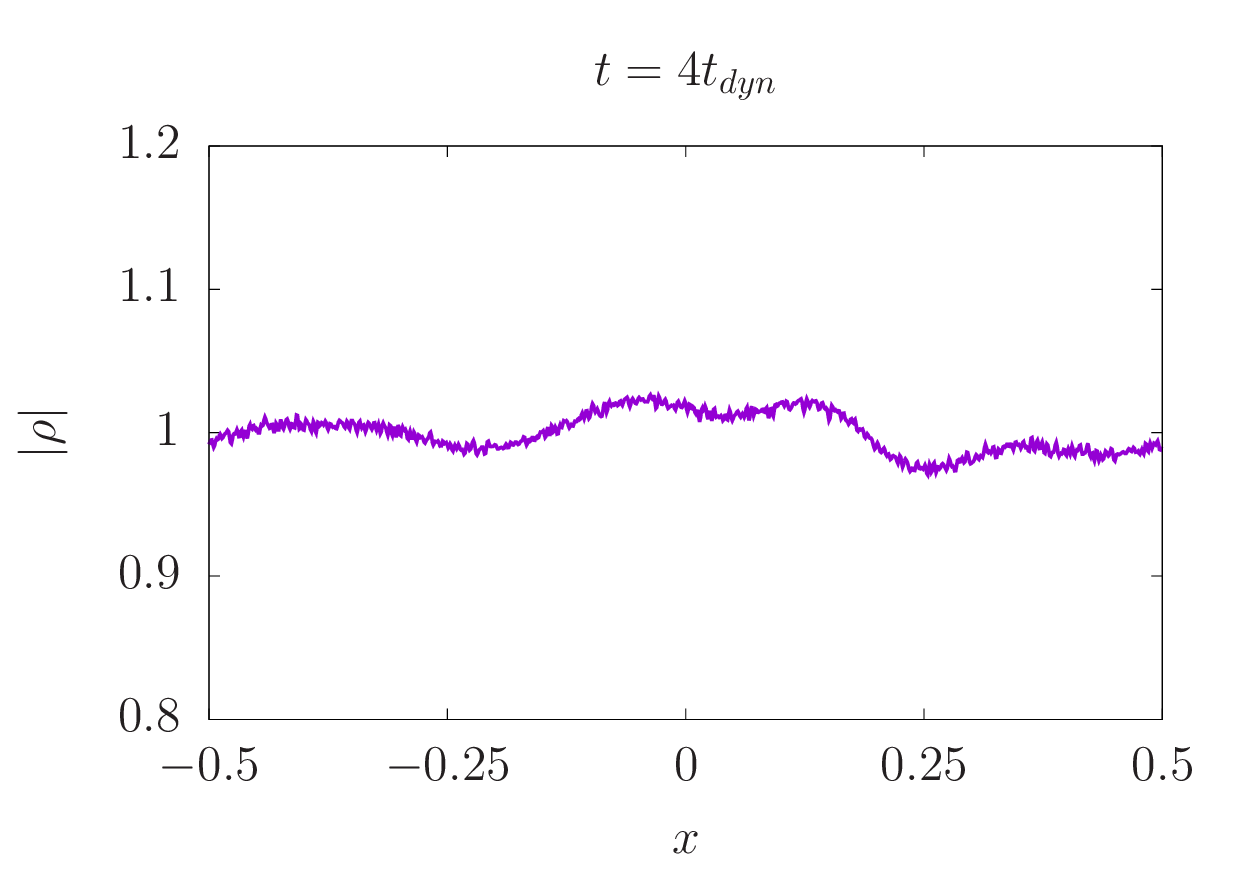}
\includegraphics[width=0.32\textwidth]{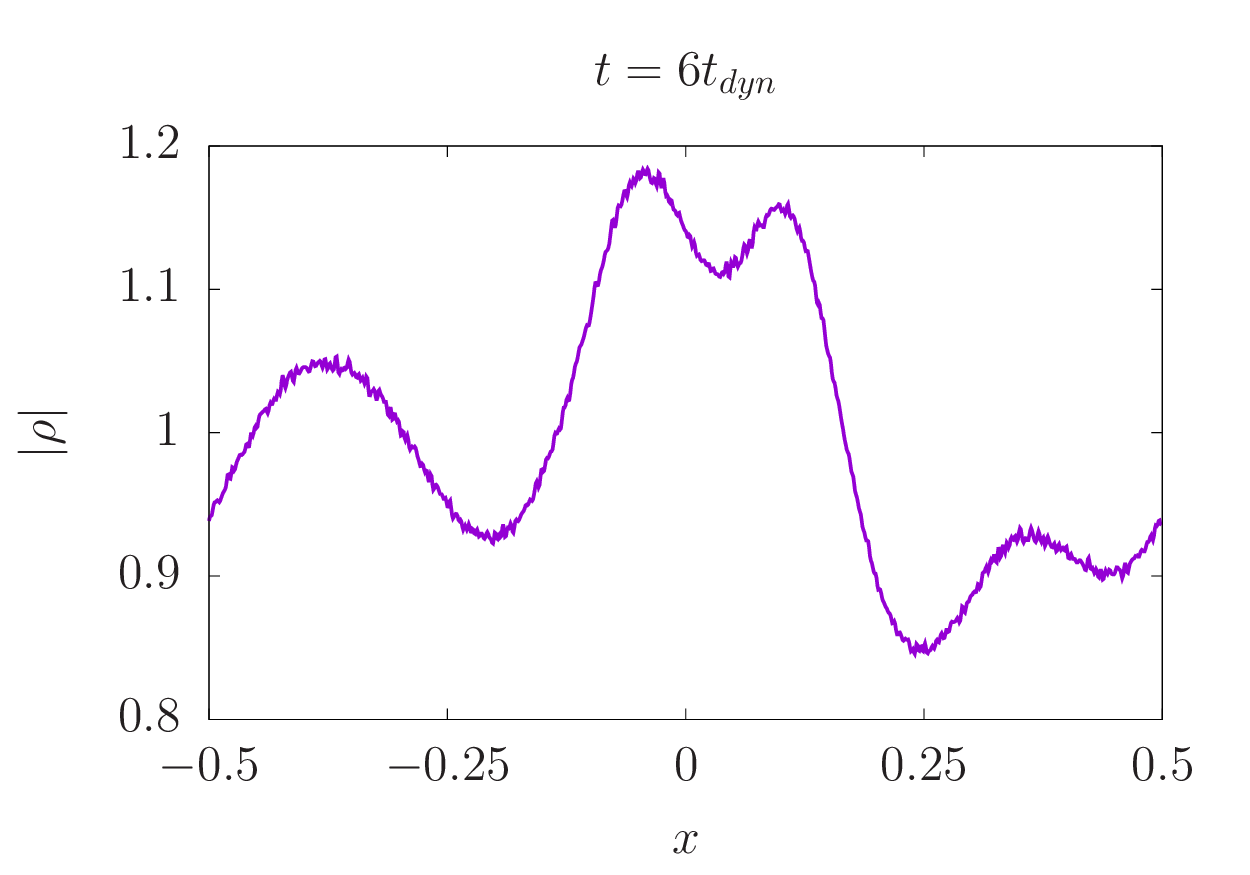}

	\caption{Snapshots of the modulus squared of the solution of the 1D SN equation (periodical case) $\left| \psi \right|^2$.}
	\label{solsp1dc}
\end{figure}

\begin{figure}
	\resizebox{0.47\textwidth}{!}{\includegraphics{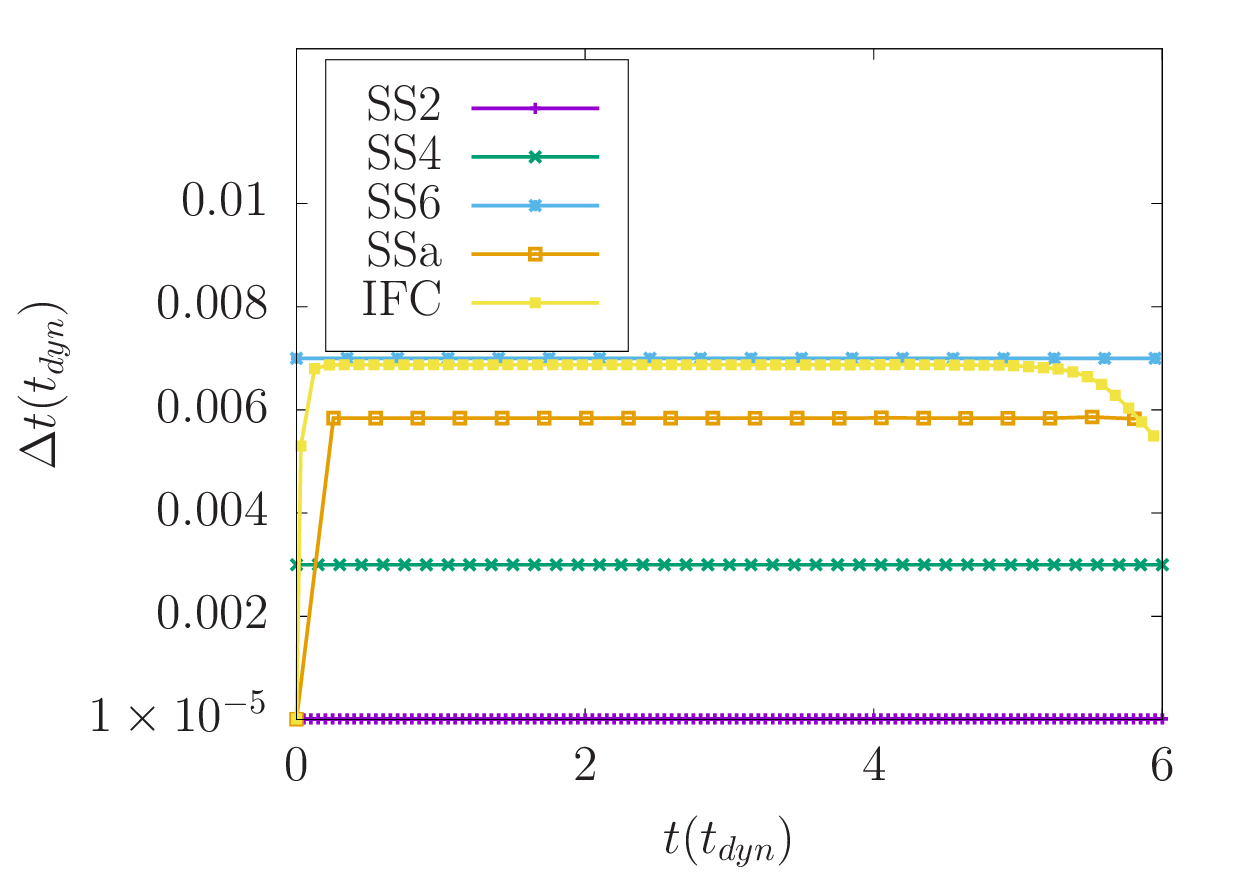}}
	\resizebox{0.47\textwidth}{!}{\includegraphics{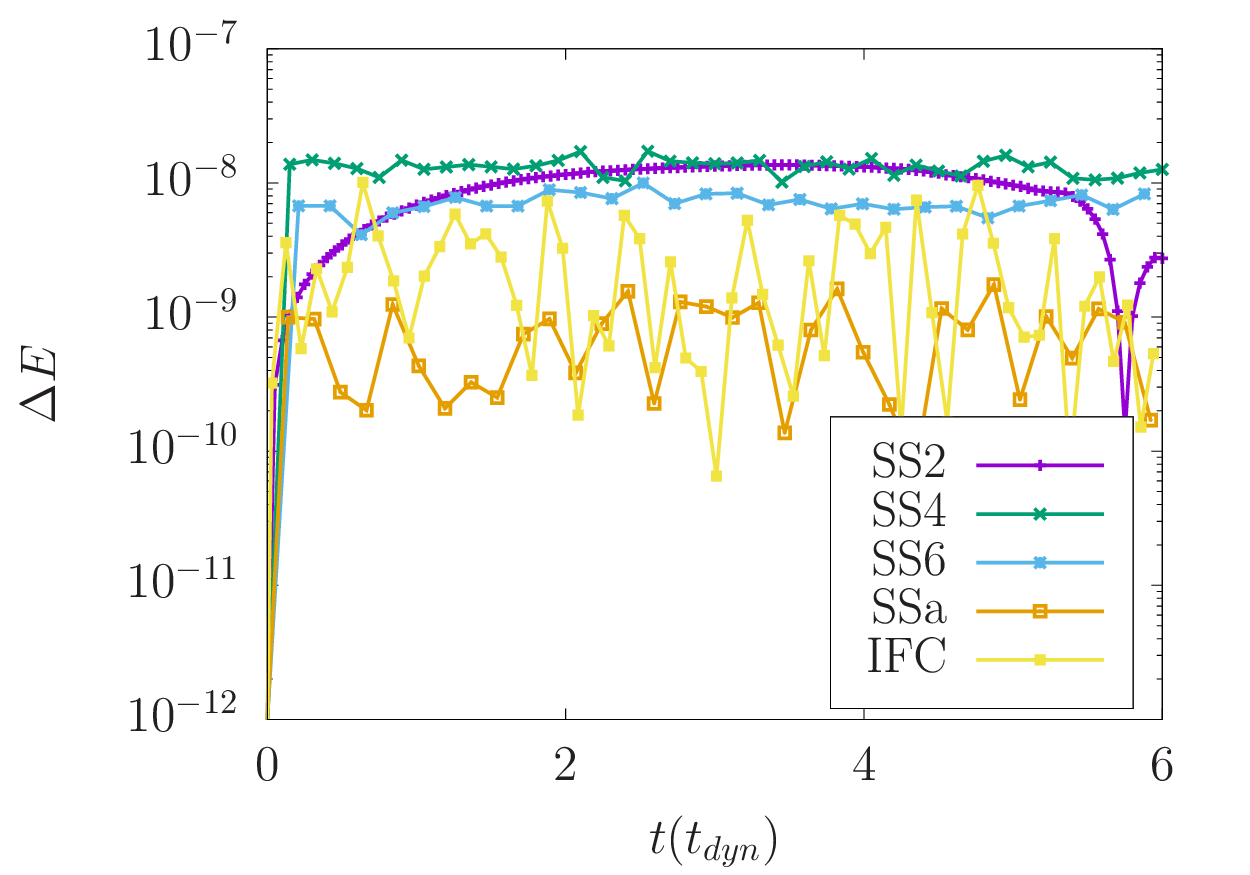}}
	\caption{Comparison for the 1D SN equation (periodical case) between the time-step (left plot) and the error on the energy conservation (right plot) for the IFC method and the Split-Step solvers.}
	\label{sp1dcosmof2}
\end{figure}

In Fig.~(\ref{solsp1dc}), the typical evolution of the system in 
the cosmological context is shown: the initial condition is spatially 
homogeneous with small fluctuations. The fluctuations grow due to gravitational 
interactions, up to be dominated by the finite size of the simulation box. 
The characteristic time of dynamics is defined as 
$t_{\mathrm{dyn}}=\left|g\right|^{-1/2}$.
In Fig.~\ref{sp1dcosmof2}, the Split-Step integrators are compared with the 
IFC. 
We observe for $t\gtrsim 5\/t_{\mathrm{dyn}}$ that the time-step decreases; 
this is due to the fact that the dynamics switches from a regime where the 
largest scales are still linear, to a regime where all the scales are nonlinear
\cite{peebles}. It indicates that the integrating factor is particularly 
efficient in the weakly nonlinear regime, which is the regime of interest 
in cosmological simulations. The Split-Step integrators (except SS2) are observed to perform in 
the same manner in the weakly non-linear and strongly non-linear regime. We observe that IFC outperforms the tested 
Split-Step integrators in the first regime, whereas, in the second one it 
becomes equally efficient compared to the split-step methods. This is 
consistent with the observation of Sect.~\ref{sect-1D-SN}: since the 
dynamics corresponds to a highly nonlinear regime, the
Split-Step method performs better than the IFC one in this case. Looking to 
Table (\ref{sp1dcosmot2}), it is clear that (for the whole simulation of 
this system) the IFC is the most efficient integration method. 
\begin{table}
	\centering
	\begin{tabular}{ |c|c|c|c|c|c| } 
		\hline
		Method & $\Delta t /t_{\mathrm{dyn}}$ & $\overline{\Delta E}$ & $\Delta \psi_{\mathrm{rev}}$ & $\Delta \psi_{\mathrm{ref}}$ & $T(s)$ \\
		\hline
		SS2 & $ 10^{-5}$& $1.6 \cdot 10^{-8}$ & $8.3 \cdot 10^{-9}$ & $ 8.5 \cdot 10^{-10}$ & 276.5 \\ 
		\hline
		SS4 & $3 \cdot 10^{-3}$& $1.9 \cdot 10^{-8}$& $1.3 \cdot 10^{-9}$& $ 1.5 \cdot 10^{-11}$ & 1.8 \\ 
		\hline
		SS6 & $7 \cdot 10^{-3}$& $1.1 \cdot 10^{-8}$& $1.1 \cdot 10^{-8}$ & $ 3.1 \cdot 10^{-11}$& 2.3\\ 
		\hline
		SSa, $\mathrm{tol}=10^{-9}$ &$5.8 \cdot 10^{-3}$& $ 2.3 \cdot 10^{-9}$ & $3.3 \cdot 10^{-9}$& $1.2 \cdot 10^{-11}$& 4.3\\ 
		\hline
		IFC, $\mathrm{tol}=10^{-12}$ &$6.7 \cdot 10^{-3}$& $ 1.1 \cdot 10^{-8}$ & $1.8 \cdot 10^{-10}$ & $6.7 \cdot 10^{-11}$ & 1.3\\ 
		\hline
	\end{tabular}
	\caption{Comparison for the 1D SN equation (periodical case) between the IFC method and the Split-Step solvers. $T$ is the total time required to run each simulation, measured in seconds. The $\Delta t$ for adaptive algorithms is the averaged one.}
	\label{sp1dcosmot2}
\end{table}

\subsection{2D nonlinear Schr{\"o}dinger equation}

In the 2D NLS case, one has
\begin{equation}\label{nls-2D}
\ui\,\partial_t\,\psi\ +\ \half\,\nabla^{2}\,\psi\ 
-\ g\,\left| \psi \right|^2\,\psi\ =\ 0.
\end{equation}
The dynamics of \eqref{nls-2D} presents a finite time singularity: it can be 
proven \cite{sulem2007nonlinear} that there exists a finite time when the norm 
of the solution or of one of its derivatives becomes infinity. This happens 
whenever the initial condition $\psi_0$ satisfies 
$E_g=\half\int \ud \boldsymbol{r}\, \psi_0 \left(g\left|\psi_0\right|^2-
\nabla^2 \right) \psi^*_0 < 0$.  The initial condition is taken as
$\psi(\boldsymbol{r}, t=0)=\ue^{-r^2/2}/\sqrt{\pi}$ and we study the cases 
$g=-1$ and $g=-6$ with respective initial energies $E_{g=-1} 
\approx 0.42$ and $E_{g=-6} \approx 0.02$. Thus, the latter is associated 
with an initial energy closer to the singular regime than the former.  
The simulation is run in a box of side $L=80$, discretized into $N=1024
\times1024$ points in the $g=-1$ case, while for $g=-6$ we set $L=120$ and  
$N=4096\times4096$. The characteristic time of dynamics is defined as 
$t_{\mathrm{dyn}}=\left|g\right|^{-1/2}$. 

\begin{figure}
	\minipage{0.47\textwidth}
	\resizebox{\textwidth}{!}{\includegraphics{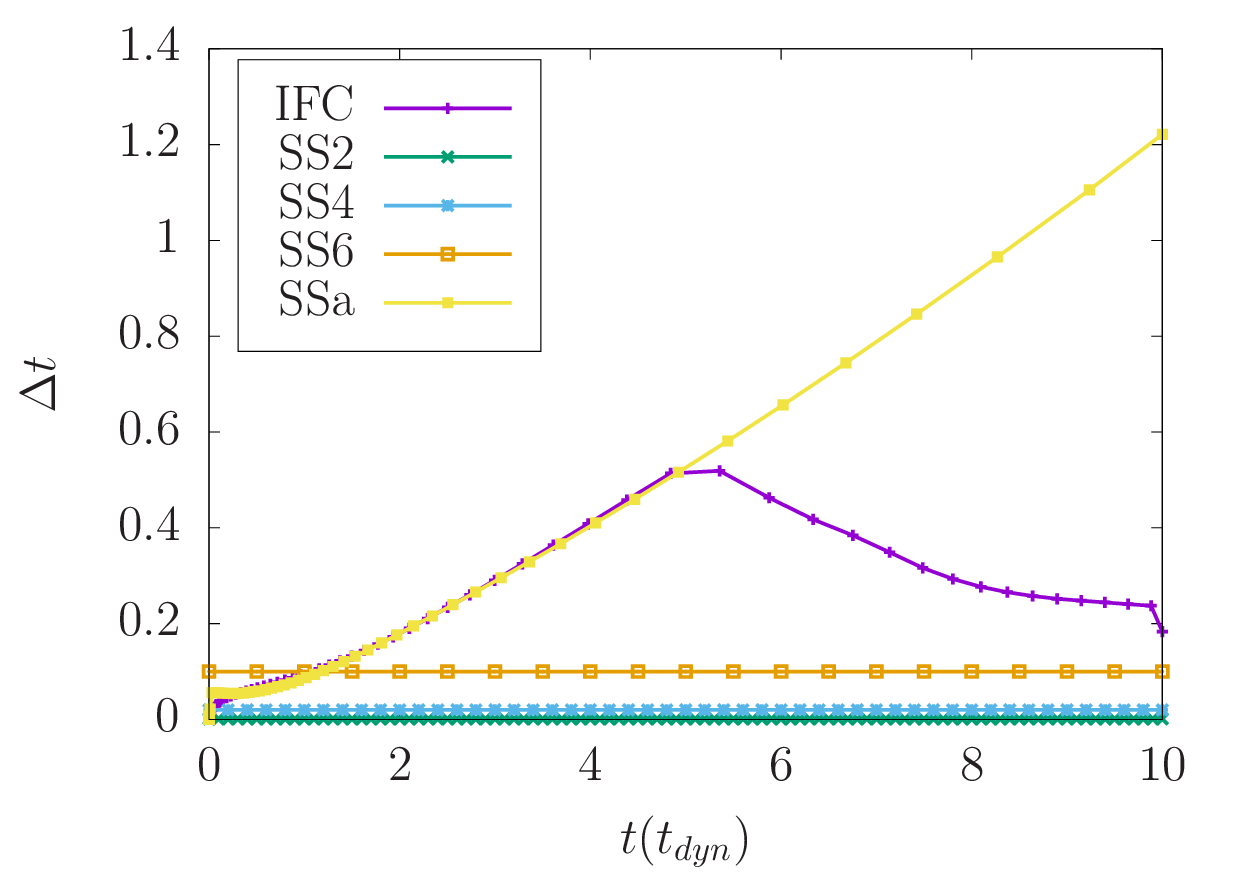}}
	\endminipage\hfill
	\minipage{0.47\textwidth}
	\resizebox{\textwidth}{!}{\includegraphics{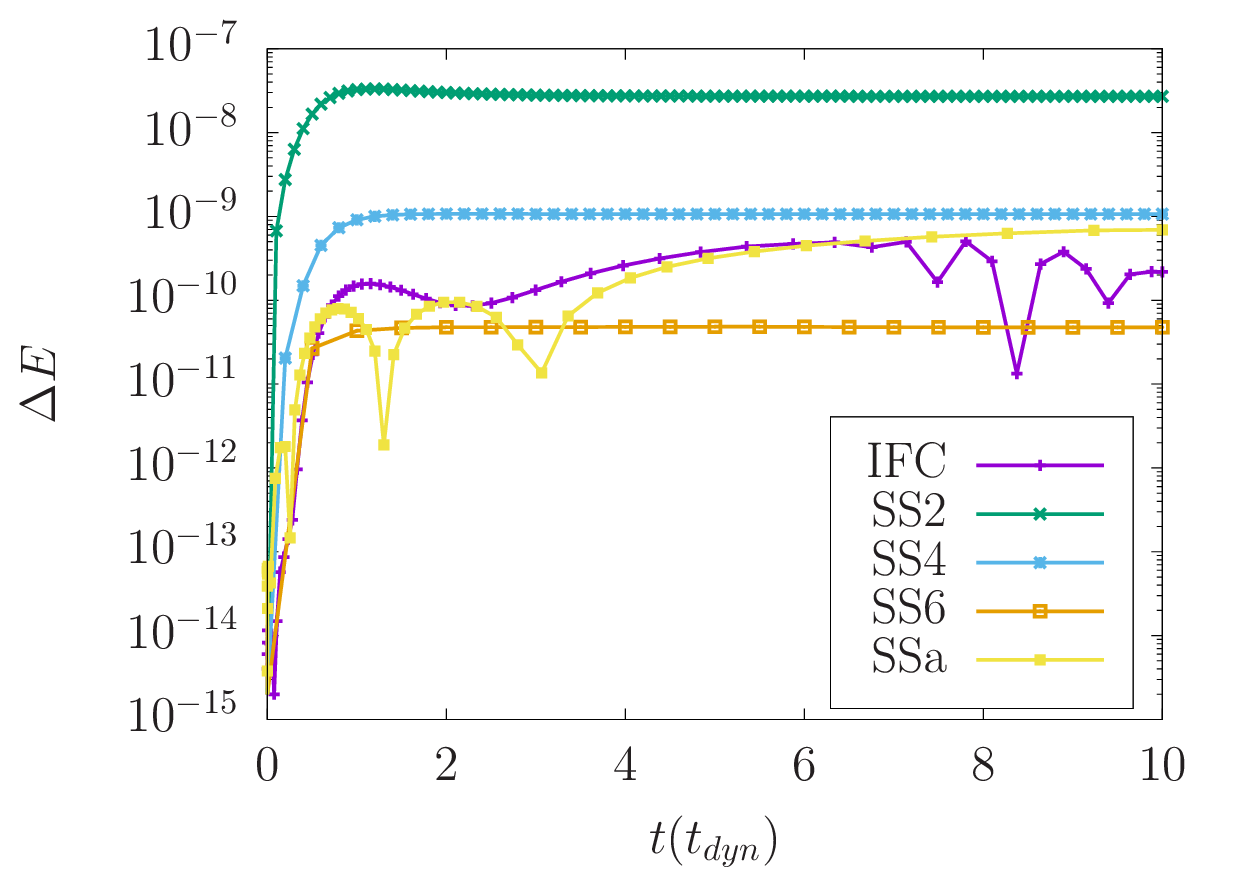}}
	\endminipage
	
	\minipage{0.47\textwidth}
	\resizebox{\textwidth}{!}{\includegraphics{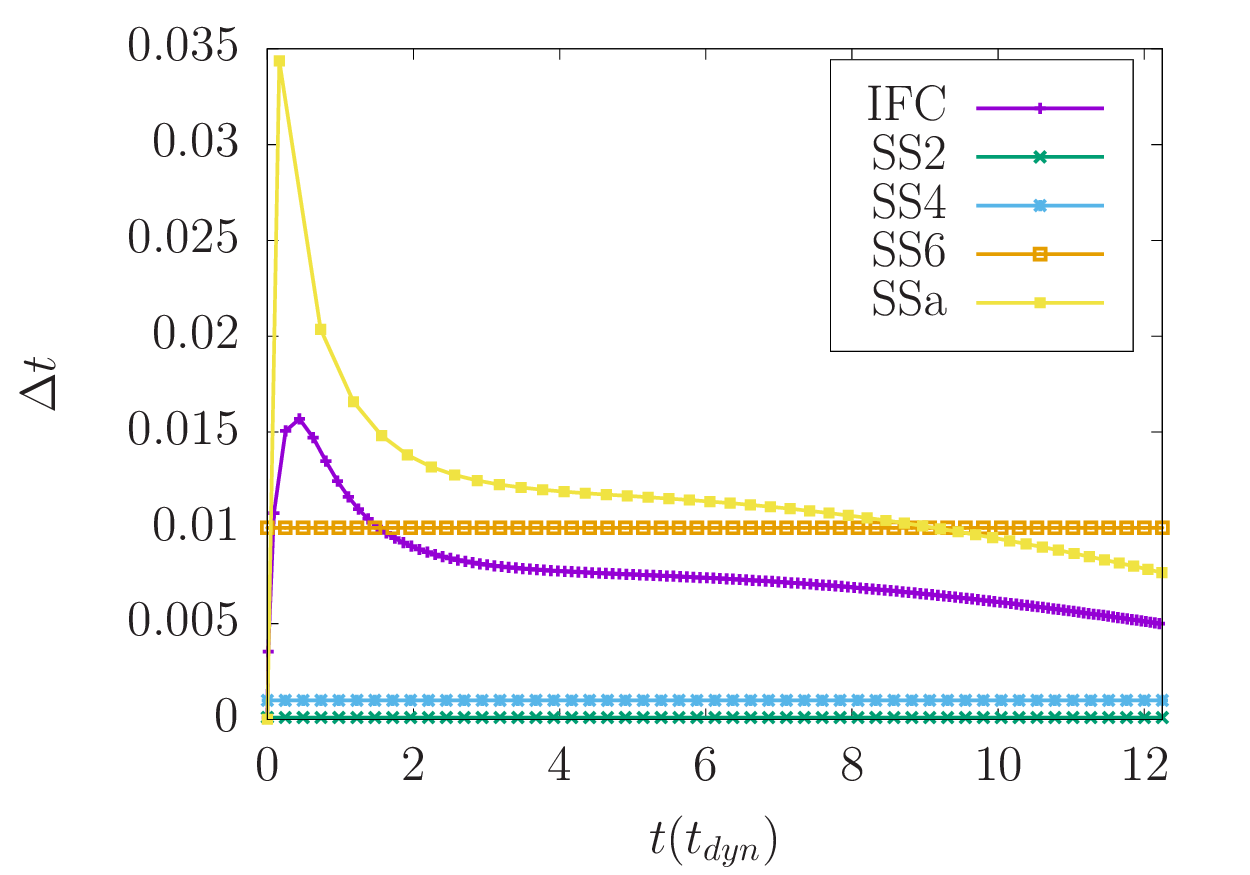}}
	\endminipage\hfill
	\minipage{0.47\textwidth}
	\resizebox{\textwidth}{!}{\includegraphics{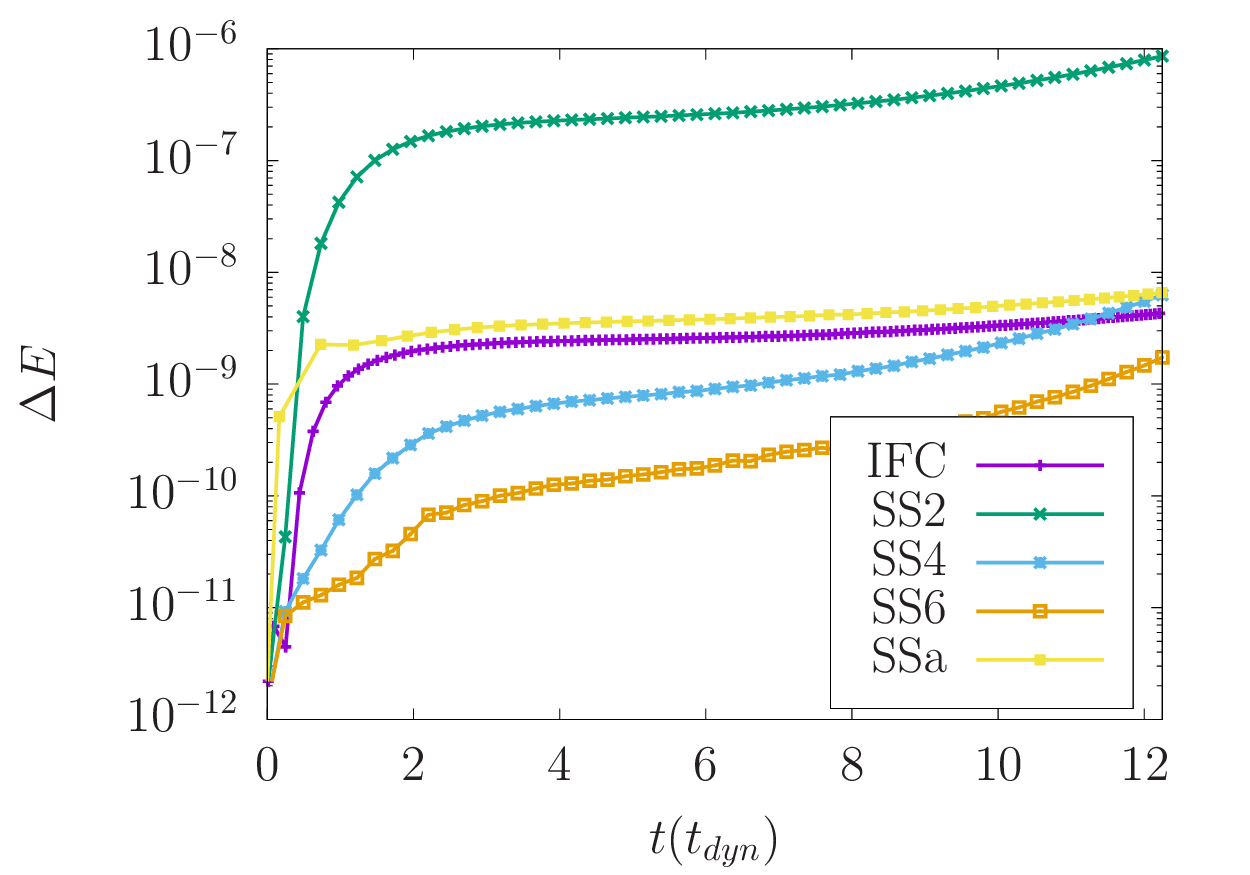}}
	\endminipage
	\caption{Comparison for the 2D NLS equation between the time-step (left panel) and the error on the energy conservation (right panel) for the IFC method and the Split-Step solvers. The top row corresponds to the $g=-1$ case while the bottom row to $g=-6$.  The last time-step is chosen in such a way that the final time is the same in all simulations, in order to ensure that $\Delta \psi_{\mathrm{rev}}$ and $\Delta \psi_{\mathrm{ref}}$ are evaluated properly, which explain the step which appears for the last time point in $\Delta t$.}
	\label{nls2df2}
\end{figure}

Split-Step integrators are compared with the IFC, looking at the energy 
conservation error, the error on the solution and the total time needed 
to run each simulation. The results are illustrated in Fig.~(\ref{nls2df2}) 
and table (\ref{nls2dt2}). 
\begin{table}
	\centering
	\begin{tabular}{ |c|c|c|c|c|c|c| } 
		\hline
		g& Method & $\Delta t/t_{\mathrm{dyn}}$ & $\overline{\Delta E}$ & $\Delta \psi_{\mathrm{rev}}$ & $\Delta \psi_{\mathrm{ref}}$ & $T(s)$ \\
		\hline
		\multirow{3}{3em}{-1}&SS2 & $ 10^{-3}$& $3.3 \cdot 10^{-8}$ & $ 4.5 \cdot 10^{-8}$ & $5.0 \cdot 10^{-10}$ &6939\\ 
		&SS4 & $2 \cdot 10^{-2}$&  $1.1 \cdot 10^{-9}$ & $3.5 \cdot 10^{-13}$ &$7.0 \cdot 10^{-11}$  &830\\ 
		&SS6 & $ 10^{-1}$ & $ 4.5 \cdot 10^{-11}$ &$ 2.5 \cdot 10^{-14}$ & $ 6.2 \cdot 10^{-11}$ & 445\\
		&SSa & $2.3 \cdot 10^{-1}$& $7.0 \cdot 10^{-10}$ & $3.0 \cdot 10^{-10}$& $2.5 \cdot 10^{-11}$ & 267 \\ 
		&IFC &$ 1.7 \cdot 10^{-1}$& $ 4.5 \cdot 10^{-10}$ &$3.0 \cdot 10^{-10}$ & $1.6 \cdot 10^{-11}$ &169\\ 
		\hline
		\multirow{3}{3em}{-6}&SS2 & $ 2.5 \cdot 10^{-4}$& $ 8.8 \cdot 10^{-7}$ & $  1.4 \cdot 10^{-11}$ & $7.9  \cdot 10^{-7}$ & 405012\\ 
		&SS4 & $ 2.5 \cdot 10^{-3}$&  $6.3 \cdot 10^{-9}$ & $2.5 \cdot 10^{-12}$ & $3.7 \cdot 10^{-9}$  &82891\\ 
		&SS6 & $ 2.5 \cdot 10^{-2}$& $1.7 \cdot 10^{-9}$ & $ 1.6 \cdot 10^{-12}$ & $1.8 \cdot 10^{-9}$ & 24453\\
		&SSa & $2.8 \cdot 10^{-2}$& $6.5 \cdot 10^{-9}$ & $1.4 \cdot 10^{-10}$&  $5.4 \cdot 10^{-9}$ & 29117\\ 
		&IFC & $ 1.8 \cdot 10^{-2}$& $ 2.3\cdot 10^{-9}$  &  $ 3.5\cdot 10^{-9}$& $ 9.4\cdot 10^{-9}$ &22843\\ 
		\hline
	\end{tabular}
	\caption{Comparisons for the 2D NLS equation between different methods for the Dormand and Prince integrator. The $\Delta t$ for adaptive algorithms is the averaged one. $T$ is the total time required to run each simulation, measured in seconds. The tolerances of the integrator SS4(3) is $\mathrm{tol}=10^{-6}$ and $\mathrm{tol}=10^{-10}$ for the IFC.}
	\label{nls2dt2}
\end{table}
The gain factor between splitting algorithms and the IFC method depends on 
the value of $g$. However, in both cases, the optimized integrating factor 
proved to be more efficient.

\subsection{2D Schr{\"o}dinger-Newton equation}

In the 2D SN case, one has
\begin{equation}
\begin{gathered}
\label{sp2d}
\ui\,\partial_t\,\psi\ +\ \half\,\nabla^{2}\psi\ 
-\ V\, \psi\ =\ 0, \qquad \nabla^2 V\ =\ g \left| \psi \right |^{2}.
\end{gathered}
\end{equation}
Similarly to the one dimensional case, we use a Gaussian initial conditions 
$\psi(\boldsymbol{r}, t=0)=\ue^{-r^2/2}/\sqrt{\pi}$ and two values of the 
coupling constant, $g=10$ and $g=500$. The former corresponds to a system in 
the quantum regime and the latter is closer to the semi-classical one. The 
potential $V$, as in the 1D case, is calculated using the Hockney method 
\cite{hockney2021computer}.  The simulation is run in a box of side $L=40$, 
discretized into $N=1024^2$ points in the $g=10$ case, while for $g=500$ we 
set $L=20$ and $N=1024\times1024$, the characteristic time of dynamics is 
defined as $t_{\mathrm{dyn}}=\left|g\right|^{-1/2}$.

\begin{figure}

	\minipage{0.47\textwidth}
	\resizebox{\textwidth}{!}{\includegraphics{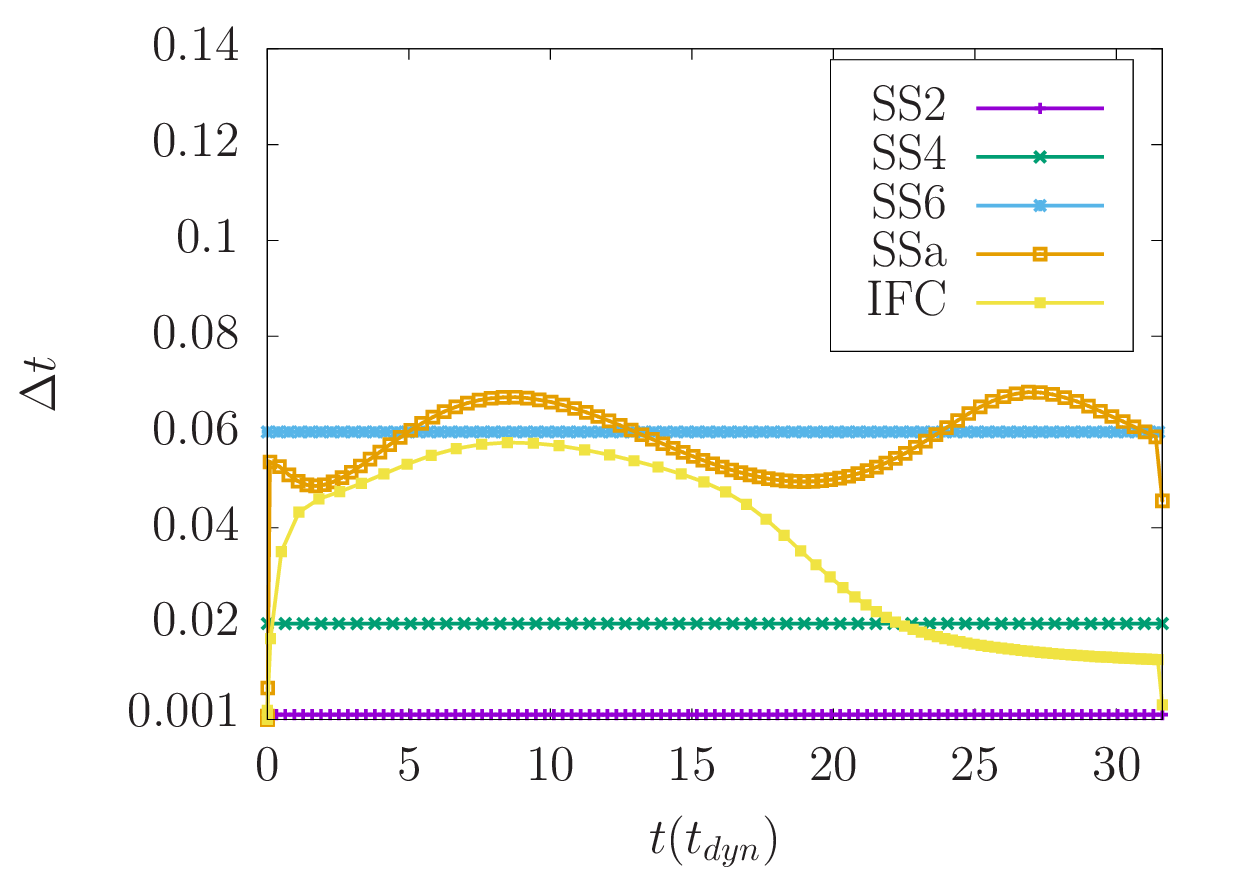}}
	\endminipage\hfill
	\minipage{0.47\textwidth}
	\resizebox{\textwidth}{!}{\includegraphics{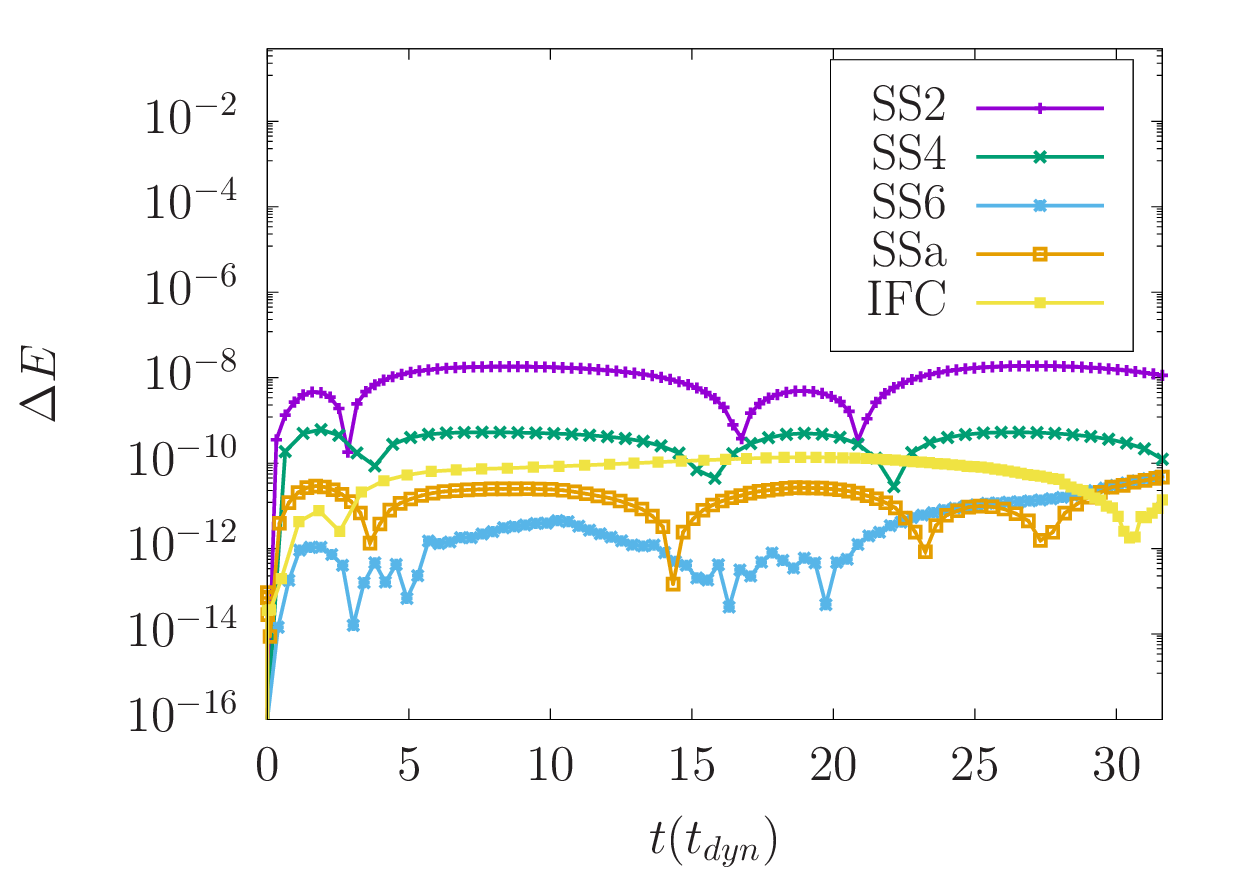}}
	\endminipage
	
	\minipage{0.47\textwidth}
	\resizebox{\textwidth}{!}{\includegraphics{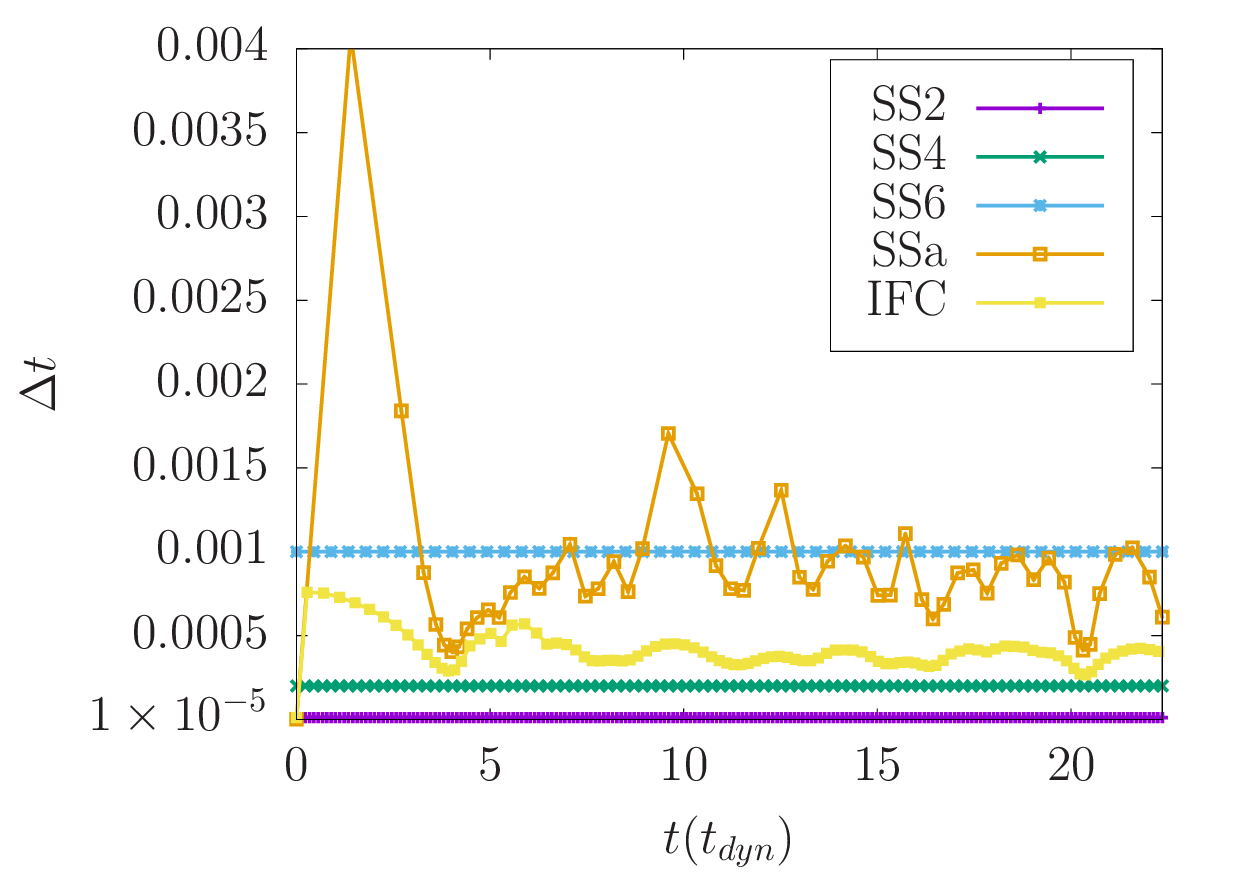}}
	\endminipage\hfill
	\minipage{0.47\textwidth}
	\resizebox{\textwidth}{!}{\includegraphics{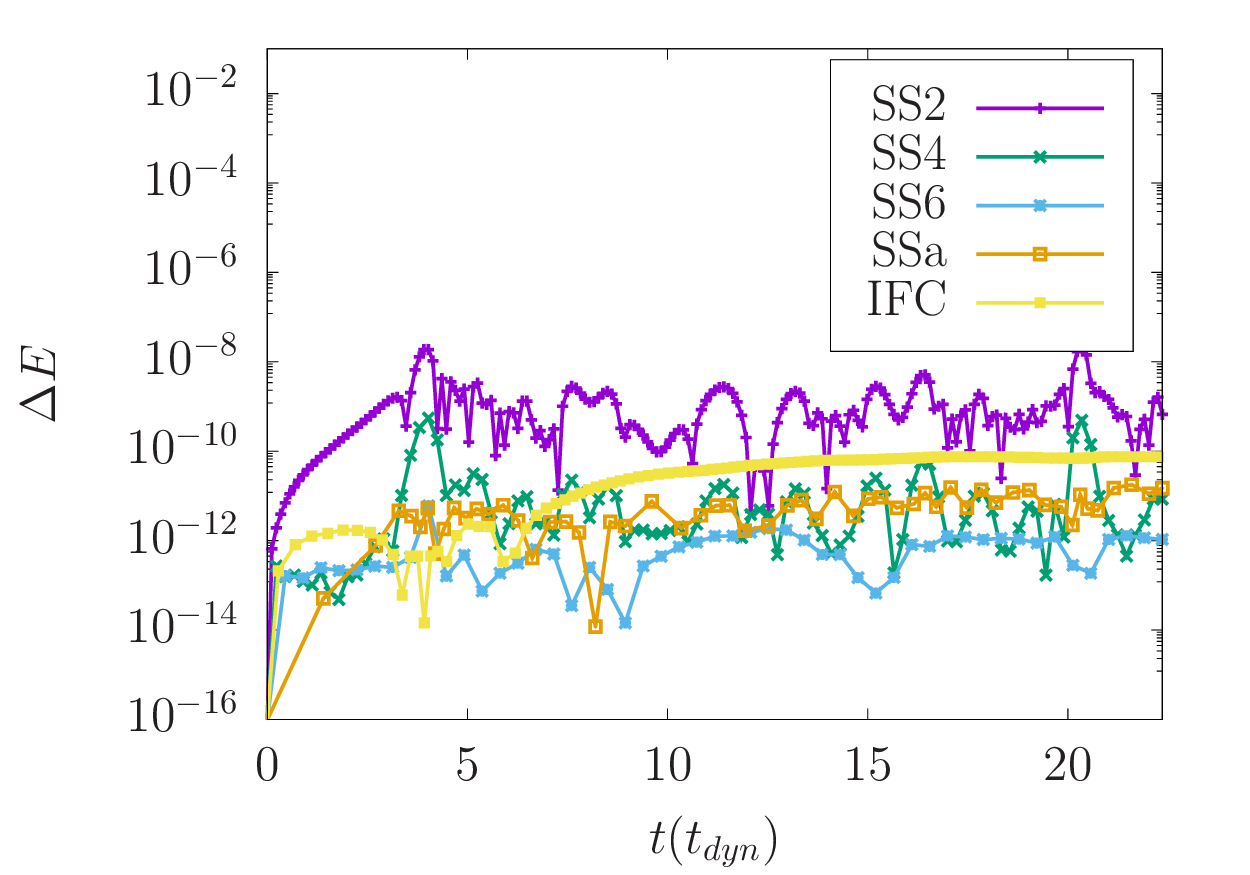}}
	\endminipage
	
	\caption{Comparison for the 2D SN equation between the time-step and the error on the energy conservation with the IFC method and the Split-Step solvers for both the cases $g=10$ (upper plots) and $g=500$ (lower plots).}
	\label{sp2df2}
\end{figure}

In table (\ref{sp2dt2}) and Fig.~(\ref{sp2df2}), we compare the Split-Step 
and the IFC integrators, looking at the energy conservation error, the error 
on the solution and the total time needed to run each simulation.

\begin{table}
	\centering
	\begin{tabular}{ |c|c|c|c|c|c|c| } 
		\hline
		g& Method & $\Delta t/t_{\mathrm{dyn}}$ & $\overline{\Delta E}$ & $\Delta \psi_{\mathrm{rev}}$ & $\Delta \psi_{\mathrm{ref}}$ & $T(s)$ \\
		\hline
		\multirow{5}{3em}{10}&SS2 & $ 3.2 \cdot 10^{-3}$& $1.8 \cdot 10^{-8}$ & $ 3.5 \cdot 10^{-12}$ & $ 1.0 \cdot 10^{-7}$ & 9070\\ 
		&SS4 & $6.3 \cdot 10^{-2}$&  $ 6.0 \cdot 10^{-10}$ & $1.4 \cdot 10^{-11}$& $ 3.8 \cdot 10^{-8}$ & 959\\ 
		&SS6 & $ 1.9 \cdot 10^{-1}$ & $ 4.6 \cdot 10^{-11}$ & $ 1.6 \cdot 10^{-7}$ & $ 4.1 \cdot 10^{-8}$ & 932 \\
		&SSa &$ 1.7 \cdot 10^{-1}$&$ 4.7 \cdot 10^{-11}$  & $ 1.4 \cdot 10^{-7}$ & $ 3.9 \cdot 10^{-8}$ &1174 \\ 
		&IFC&$ 8.4 \cdot 10^{-2}$& $ 1.3 \cdot 10^{-10}$ &$5.8 \cdot 10^{-8}$ & $ 2.9 \cdot 10^{-8}$ & 1172\\ 
		\hline
		\multirow{5}{3em}{500}&SS2 & $ 2.2 \cdot 10^{-4}$& $2.0 \cdot 10^{-8}$ & $3.1 \cdot 10^{-10}$ & $9.1 \cdot 10^{-8}$ & 84030 \\ 
		&SS4 & $4.5 \cdot 10^{-3}$&  $5.6 \cdot 10^{-10}$& $2.8 \cdot 10^{-11}$ & $ 7.0 \cdot 10^{-9}$ & 8401\\ 
		&SS6 & $2.2 \cdot 10^{-2}$ & $6.1 \cdot 10^{-12}$ & $9.1 \cdot 10^{-11}$ & $4.1 \cdot 10^{-9}$ & 5143 \\
		&SSa &$2.1 \cdot 10^{-2}$& $1.9 \cdot 10^{-11}$ &$2.0 \cdot 10^{-11}$&$3.9 \cdot 10^{-9}$ & 6676\\ 
		&IFC &$9.1 \cdot 10^{-3}$& $7.5 \cdot 10^{-11}$ & $3.0 \cdot 10^{-10}$ & $ 1.6 \cdot 10^{-9}$ & 6637\\ 
		\hline
	\end{tabular}
	\caption{Comparison for the 2D SN equation between the IFC method and the Split-Step solvers. The $\Delta t$ for adaptive algorithms is the averaged one. The tolerance for the SSa algorithm is $\mathrm{tol}=10^{-6}$ and $\mathrm{tol}=10^{-7}$ for $g=10$ and $g=500$ respectively, and for the IFC algorithm $\mathrm{tol}=10^{-10}$ and $\mathrm{tol}=10^{-12}$ for $g=10$ and $g=500$ respectively.}
	\label{sp2dt2}
\end{table}

For the 2D Schr{\"o}dinger--Newton equation, adaptive splitting algorithms 
proved to be as efficient as the IFC. Similarly to the one-dimensional case, 
also here the SS6 split-step algorithm with constant time-step resulted to be 
the fastest among the ones we tested. This is due to the same reasons mentioned 
in section \ref{SN1Dsect}. Note that, here, the performance gap 
between the integrating factor and splitting algorithms is smaller than in 
the one-dimensional case. Indeed, as the dynamics gets more complicated and 
the number of spatial dimensions increase, algorithms with adaptive 
time-step shall always be preferred. 

\subsubsection{Periodical case}

\begin{figure}
	
	\minipage{0.32\textwidth}
	\resizebox{\textwidth}{!}{\includegraphics{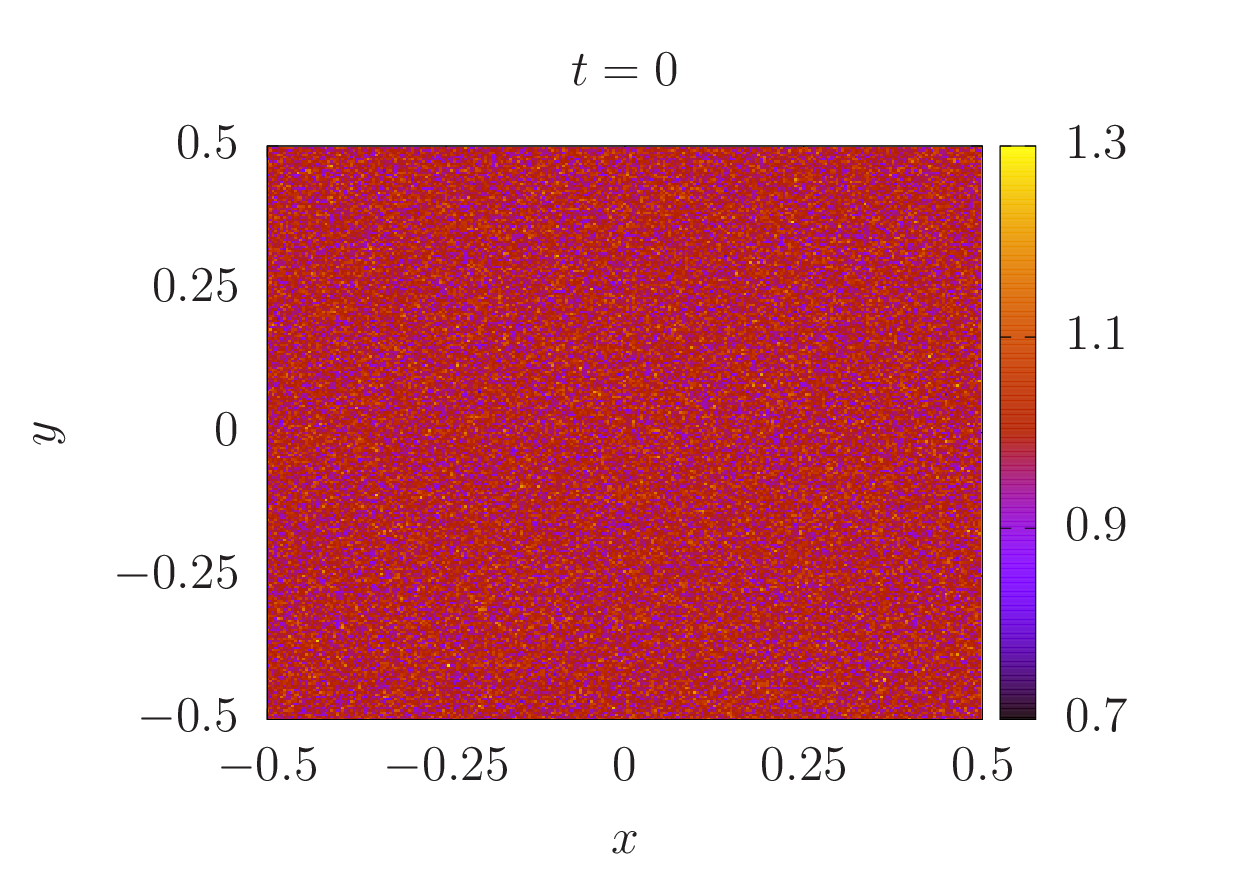}}
	\endminipage\hfill
	\minipage{0.32\textwidth}
	\resizebox{\textwidth}{!}{\includegraphics{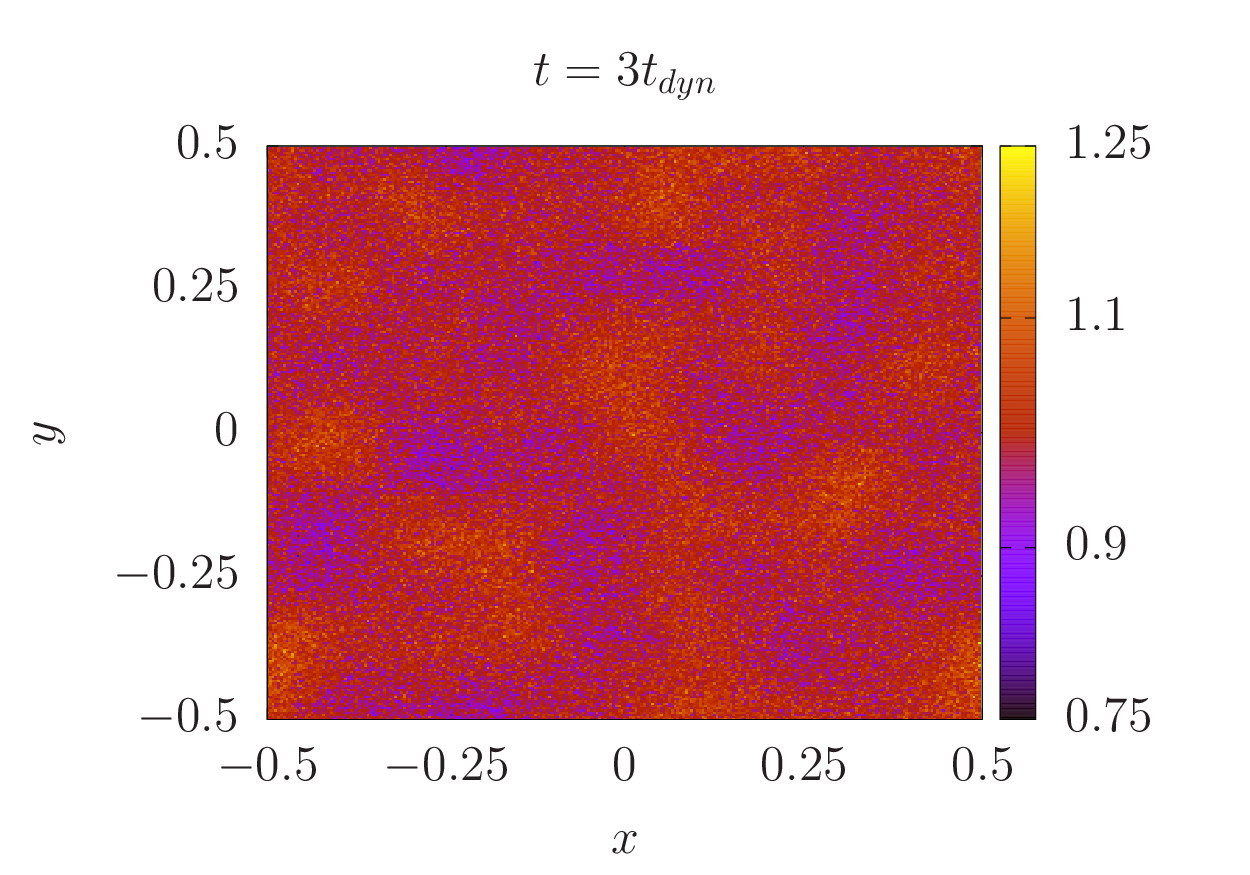}}
	\endminipage\hfill
	\minipage{0.32\textwidth}%
	\resizebox{\textwidth}{!}{\includegraphics{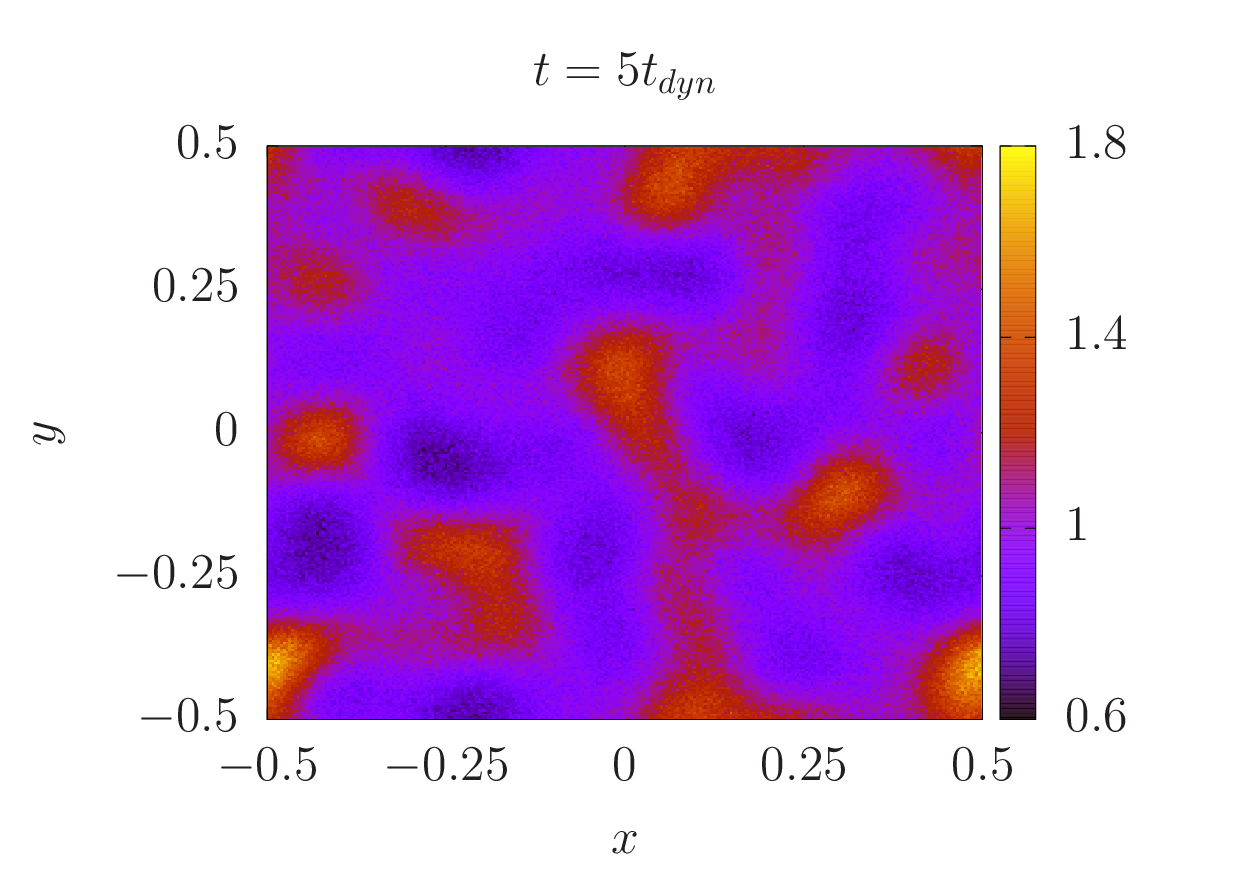}}
	\endminipage
	
	\minipage{0.32\textwidth}
	\resizebox{\textwidth}{!}{\includegraphics{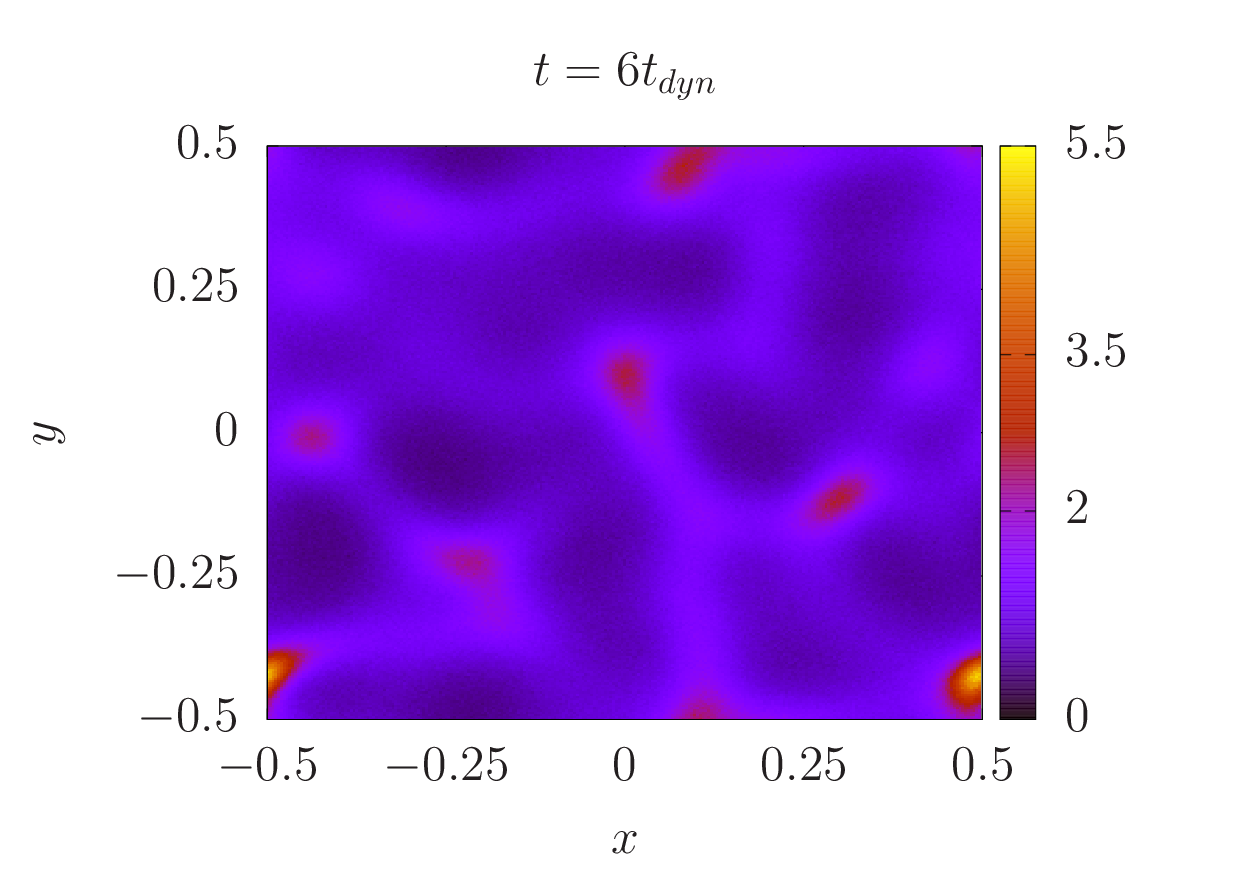}}
	\endminipage\hfill
	\minipage{0.32\textwidth}
	\resizebox{\textwidth}{!}{\includegraphics{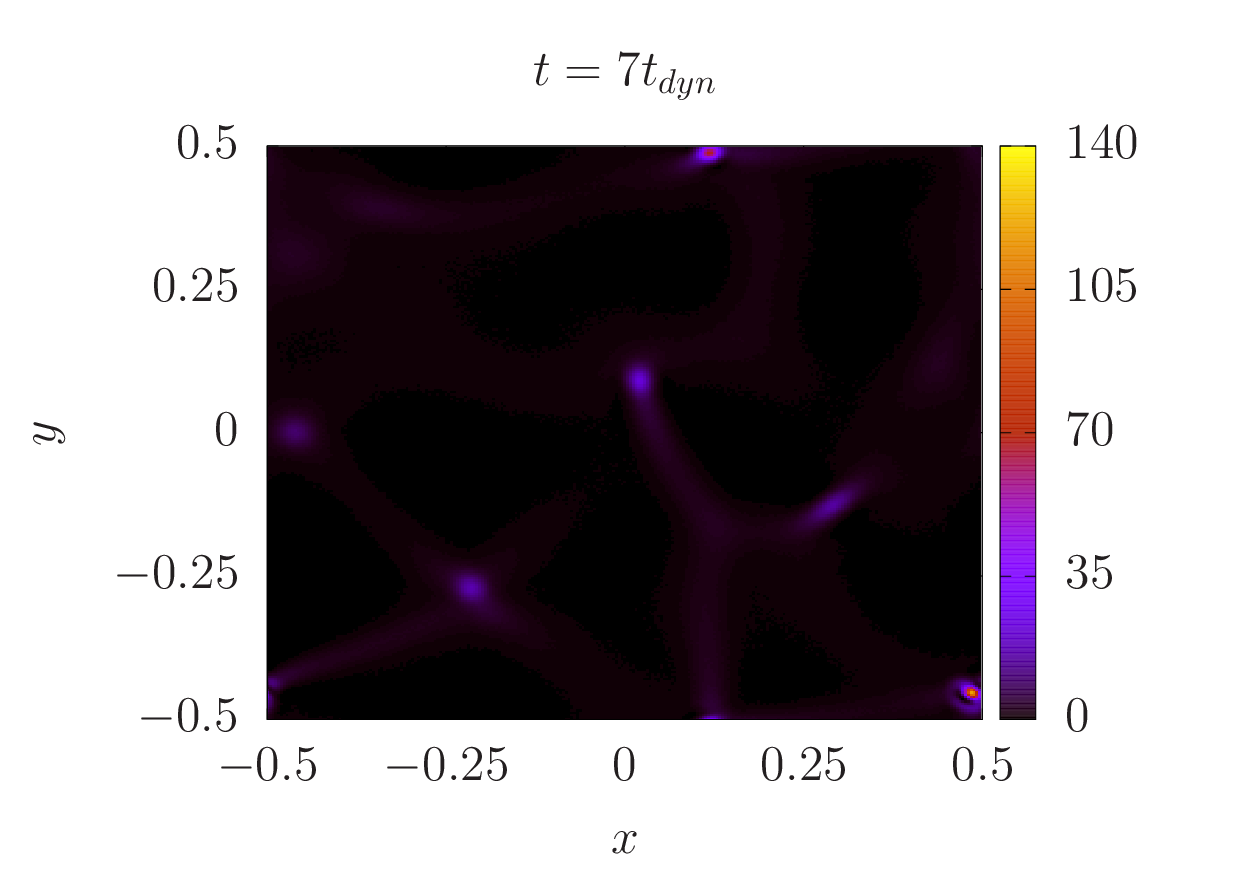}}
	\endminipage\hfill
	\minipage{0.32\textwidth}%
	\resizebox{\textwidth}{!}{\includegraphics{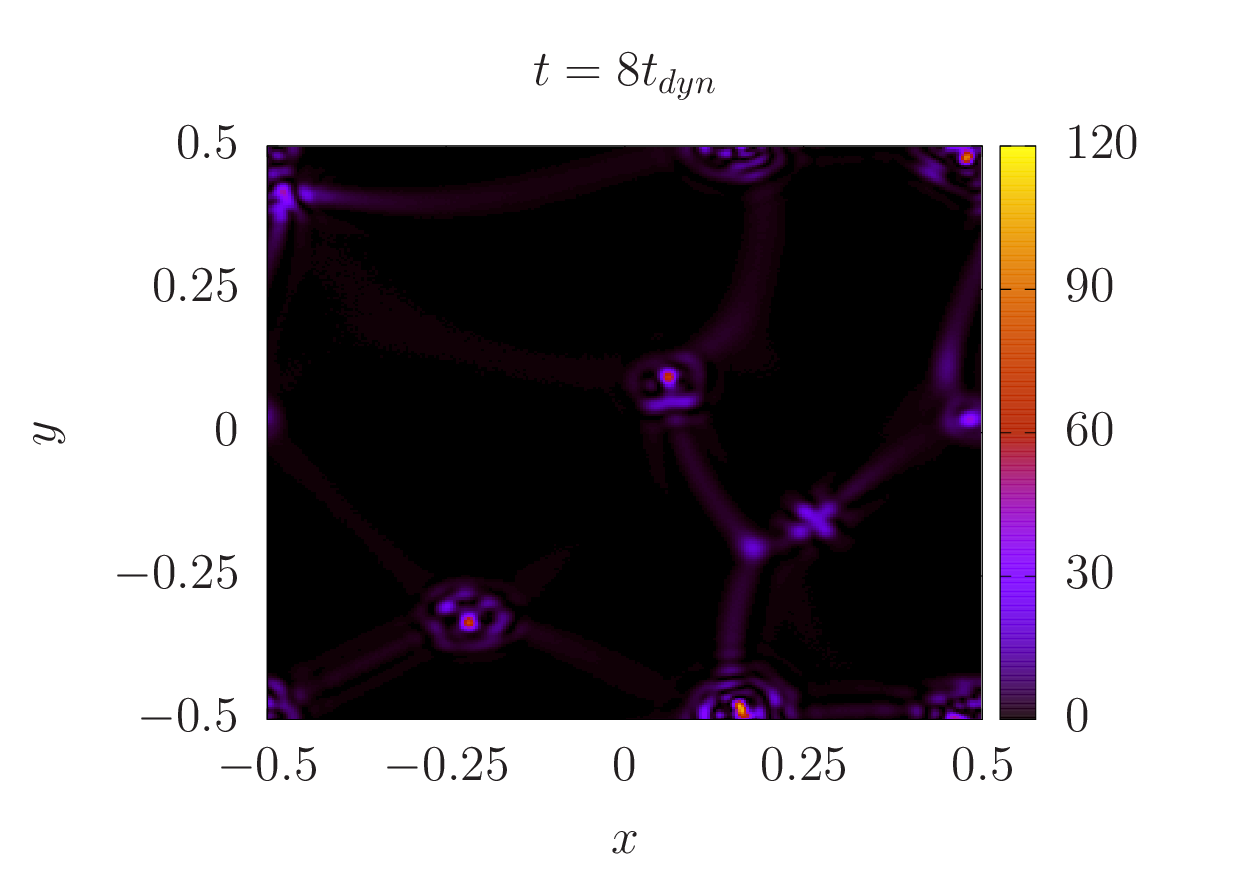}}
	\endminipage

	\caption{Snapshots of the modulus squared of the solution of the 2D SN equation (periodical case) $\left| \psi \right|^2$.}
	\label{solsp2dc}
\end{figure}

For the 2D periodical case, we run simulations in a box of side $L=1$ with 
$N=1024\times1024$, using again a constant power spectrum as initial 
condition with $g=10^6$, $\rho_0=1$ and a zero initial velocity field. In 
Fig.~(\ref{solsp2dc}) some snapshots of the modulus squared of the solution 
are shown, expressing time in units of $t_{\mathrm{dyn}}=1/\sqrt{g}$.   

\begin{figure}
	\resizebox{0.47\textwidth}{!}{\includegraphics{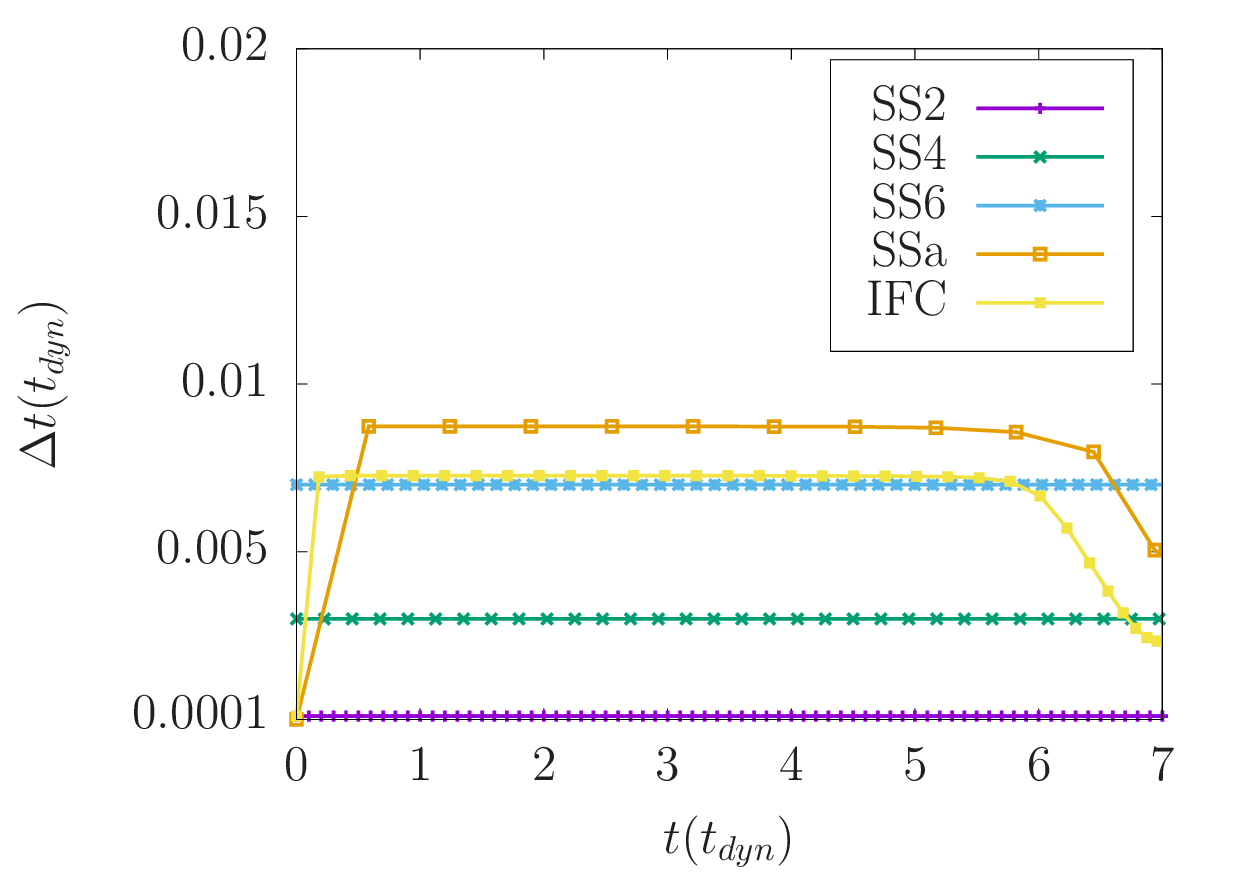}}
	\resizebox{0.47\textwidth}{!}{\includegraphics{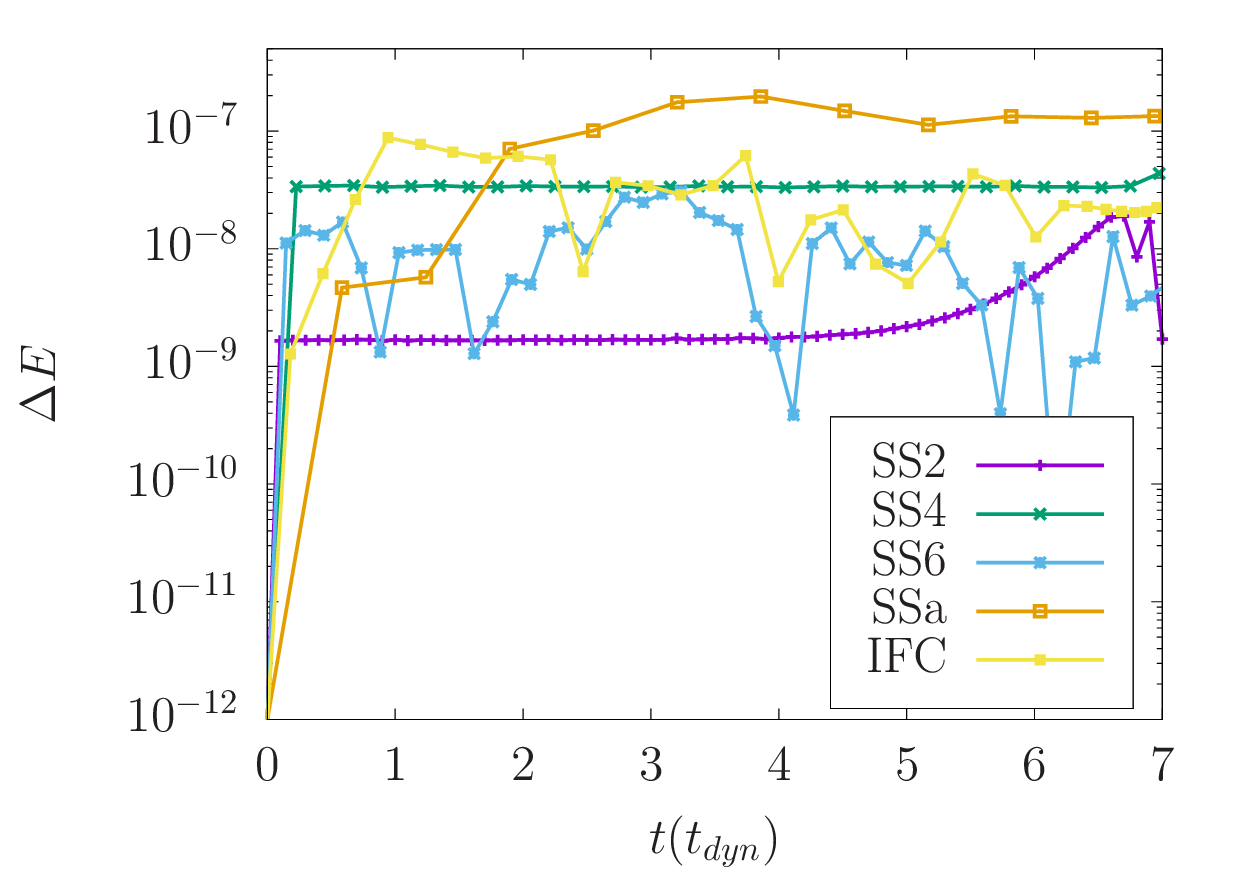}}
	\caption{Comparison for the 2D SN equation (periodical case) between the time-step (left plot) and the error on the energy conservation (right plot) for the IFC method and the Split-Step solvers.}
	\label{sp2dcosmof2}
\end{figure}

In table (\ref{sp2dcosmot2}) and Fig.~(\ref{sp2dcosmof2}), we compare the 
Split-Step and IFC  integrators, looking at the energy conservation error, 
the error on the solution and the total time needed to run each simulation. 
We obtain the same result than in one dimension, with the IFC being the most 
efficient method.
\begin{table}
	\centering
	\begin{tabular}{ |c|c|c|c|c|c| } 
		\hline
		Method & $\Delta t/t_{\mathrm{dyn}}$ & $\overline{\Delta E}$ & $\Delta \psi_{rev}$ & $\Delta \psi_{ref}$ & $T(s)$ \\
		\hline
		SS2 & $ 10^{-4}$& $1.8 \cdot 10^{-8}$ & $4.5 \cdot 10^{-9}$ & $1.0 \cdot 10^{-5}$ & 31711 \\ 
		\hline
		SS4 & $3 \cdot 10^{-3}$& $4.2 \cdot 10^{-8}$ & $3.1 \cdot 10^{-6}$ & $5.3 \cdot 10^{-6}$ & 2485 \\ 
		\hline
		SS6 & $7 \cdot 10^{-3}$ & $3.5 \cdot 10^{-8}$ & $4.5 \cdot 10^{-10}$ & $5.1 \cdot 10^{-6}$ & 3046 \\
		\hline
		SSa, $\mathrm{tol}=10^{-7}$ & $ 7.5 \cdot 10^{-3}$ & $ 2.3 \cdot 10^{-7}$  & $8.5 \cdot 10^{-7}$ & $ 1.2 \cdot 10^{-6}$ & 3362\\ 
		\hline
		IFC, $\mathrm{tol}=10^{-12}$ & $5 \cdot 10^{-3}$ & $9.5 \cdot 10^{-8}$  & $4.6 \cdot 10^{-6}$ & $8.9 \cdot 10^{-7}$ & 2232\\ 
		\hline
	\end{tabular}
	\caption{Comparison for the 2D SN equation (periodical case) between the IFC method and the Split-Step solvers. $T$ is the total time required to run each simulation, measured in seconds. The $\Delta t$ for adaptive algorithms is the averaged one.}
	\label{sp2dcosmot2}
\end{table}

\subsection{Gross--Pitaevskii--Poisson equation}

We conclude by presenting the results for the 2D Gross--Pitaevskii--Poisson 
equation, which is a combination of the NLS and SN equations
\begin{subequations}
	\begin{align}\nonumber
	&\ui\,\partial_t\,\psi\ +\ \half\,\nabla^2\,\psi\ -V\,\psi\ =\   0, \qquad
	V = V_1+V_2, \qquad
	\nabla^2 V_1\ =\ g_1 \left|\,\psi\,\right|^2, \qquad
	V_2\ =\ g_2 \left|\,\psi\,\right|^2.
	\label{gpp2d}
	\end{align}
\end{subequations}
Based on the results presented so far, in the case of open boundary 
conditions, one expects the split-step or the integrating factors to 
outperform one the other, depending on the values of $g_1$ and $g_2$. 
We set them to $g_1=-3$ and $g_2=100$ which are very close to the one 
typically employed when simulating the collapse of a self-gravitating 
Bose-Einstein condensate with attractive self-interaction 
\cite{chavanis2016collapse}. The numerical parameters are $N=2048\times 
2048$, $L=40$ and $t_{f}=5$ while the initial condition is a Gaussian, 
$\psi(\boldsymbol{r}, t=0)=\ue^{-{r^2}/{2}} / \sqrt{\pi}$. In table 
\ref{kerr2dt1} comparisons between the most efficient methods tested for 
the NLS (IFC method) and the SN in the non-periodical case (SS6 or SSa, depending on the parameters) are shown. For the values of $g_1$ and $g_2$ we use, the IFC method 
outperforms the split-step solvers. Moreover, we observe that for our 
particular initial condition, the smaller the ${g_1}/{g_2}$ ratio is, 
the better the IFC performs with respect to splitting methods, with a 
robust difference already appearing for ${g_1}/{g_2} \lessapprox 0.1$. This confirms that the presence of a short-range interaction term puts the integrating factor method in a clear more performing position, compared to splitting methods. 

\begin{table}
	\centering
	\begin{tabular}{ |c|c|c|c|c|c| } 
		\hline
		Method & $\Delta t$ & $\overline{\Delta E}$ & $\Delta \psi_{\mathrm{rev}}$ & $\Delta \psi_{\mathrm{ref}}$ & $T(s)$ \\
		\hline
		SS6 & $2.5 \cdot 10^{-3}$ & $1.2 \cdot 10^{-9}$  & $4.5 \cdot 10^{-11}$ & $4.5 \cdot 10^{-8}$  &  32753\\
		\hline
		SSa, $\mathrm{tol}=10^{-7}$ &$4 \cdot 10^{-3}$&  $2.4 \cdot 10^{-11}$&$2.2 \cdot 10^{-8}$&$1.4 \cdot 10^{-8}$ & 63315\\ 
		\hline
		IFC, $\mathrm{tol}=10^{-11}$ & $2 \cdot 10^{-3}$ & $7.7 \cdot 10^{-10}$ & $6.2 \cdot 10^{-8}$ & $2.3 \cdot 10^{-8}$ & 28273 \\ 
		\hline
	\end{tabular}
	\caption{Comparison for the 2D Gross--Pitaevskii--Poisson equation between the IFC and the Split-Step methods. $T$ is the total time required to run each simulation, measured in seconds.}
	\label{kerr2dt1}
\end{table}

\section{Conclusion} 
\label{sect-conclusions}

We studied the numerical resolution of the nonlinear Schrödinger (NLS) and 
the Newton--Schrödinger (SN) equations using the optimized integrating 
factor (IFC) technique. This method was compared with splitting algorithms. 
Specifically, for the integrating factor, we tested fifth-order time-adaptive 
algorithms, while, for the Split-Step family, we focused on second-, fourth- 
and sixth-order schemes with fixed time-step, and a fourth-order algorithm 
with adaptive time-step. We performed extensive tests with systems in one 
and two spatial dimensions, with open or periodic boundary conditions.

The comparisons between the results obtained in the tested cases, show that 
the IFC method can be more efficient than splitting algorithms, especially 
in the NLS equation and periodical SN equation cases. For the SN equation 
in the non-periodical case on the other hand, splitting algorithms proved 
to be more efficient, even though the optimized integrating factor provided 
competitive results in terms of both speed and accuracy. Moreover, the 
results obtained for the Gross--Pitaevskii--Poisson equation pointed out 
how the presence of a short-range interaction term puts the integrating 
factor method in a clear more performing position. 

Finally, the achieved results indicate how, among the splitting algorithms 
at fixed step, working with higher order solvers is always more efficient. 
In particular the Split-Step order 6 proved to be around 10 times faster 
compared with the lower order ones, while conserving the energy with the 
same error.

\addcontentsline{toc}{section}{References}
\bibliography{sample}

\appendix 
\appendixpage

\section{Optimized Integrating Factor}
\label{AppIFC}

The optimized version of the integrating factor is based on the property 
that, for the Schr{\"o}dinger equation, if the value of the potential $V$ is 
modified by an additive constant $\mathcal{C}$, only the phase of the 
solution $\psi$ is changed. Indeed, if $\psi$ is a solution of 
\eqref{def-schro} at a given time $t$, then $\Psi \eqdef
\psi\,\ue^{-\/\ui\,\mathcal{C}\,t}$ is a solution of
\begin{equation} \label{eq-gauge-rotation}
\ui\,\partial_t\,\Psi\ +\ \half\,\nabla^2 \Psi\ -\ \left(V + \mathcal{C} \right)\,\Psi\ =\ 0,
\end{equation}
as it can be easily verified.
Thus, adding a constant $\mathcal{C}_n$ to $V$, the solution is modified as
\begin{equation} \label{gauge-rotation}
\psi(t_n)\ \to\ \psi(t_n)\,\ue^{-\/\ui\/\varphi}, \qquad
\varphi\ \eqdef\ \sum_{n=0}^{N_{h}} \mathcal{C}_n\,h_n,
\end{equation}
where  $h_n\eqdef t_{n+1}-t_n$ is the $n$-th time-step and $N_{h}$ is 
the total number of time-steps. 

The freedom provided by the gauge condition of the potential is 
exploited to compute an optimal value of $\mathcal{C}_n$, which allows 
to choose a larger time-step compared to the $\mathcal{C}_n=0$ case and 
therefore speeding up the numerical integration. The resulting optimal 
value, $\tilde{\mathcal{C}}_n$,  which is obtained at each time step $n$ 
as the value of $\mathcal{C}_n$ minimising the $L_2$-norm of $\mathcal{N}$ 
\cite{lovisetto2022optimized}, is 
\begin{equation}\label{L2normres}
\tilde{\mathcal{C}}_n\,\eqdef\,-\left(\,\sum_{\ell=-\left[M/2\right]}^{\left[M/2\right]-1} V_\ell \,\left|\psi_\ell\right|^2\,\right)\,\left/\,\left(\,\sum_{\ell=-\left[M/2\right]}^{\left[M/2\right]-1} \left|\psi_\ell\right|^2\,\right)\right.
\end{equation}
where $\psi_{\ell}\eqdef \psi(\boldsymbol{r}_{\ell})$ and $V_{\ell}\eqdef V(\boldsymbol{r}_{\ell})$ at time $t_n$. 

\section{Split-Step pseudo-codes}
\label{app-ss}

We list below, in ALG. (\ref{SSNA}), the pseudo-codes for the Split-Step 
algorithms with fixed time-step. We consider the general case of order $N$, 
with $N \in \left\{  2, 4, 6\right\}$.

\begin{breakablealgorithm}
	\caption{: SSN, $N \in \left\{  2, 4, 6\right\}$}\label{SSNA}
	\begin{algorithmic}[1] 
		\State $t \gets t_0$
		\State $\psi \gets \psi(\boldsymbol{r},t_0)$
		\While{$t < t_{f}$}
		\State $\psi \gets \mathrm{FFT}^{-1}[\mathrm{exp}\left(-\ui \hat{K} a_1 h\right) \mathrm{FFT}[\psi]]$
		\State $\psi \gets \mathrm{exp}\left(- \ui V b_1 h \right)\psi$
		\State \vdots
		\State $\psi \gets \mathrm{FFT}^{-1}[\mathrm{exp}\left(-\ui \hat{K} a_{\frac{N}{2}} h\right) \mathrm{FFT}[\psi]]$
		\State $\psi \gets \mathrm{exp}\left(- \ui V b_{\frac{N}{2}} h\right) \psi$
		\State $\psi \gets \mathrm{FFT}^{-1}[\mathrm{exp}\left(-\ui \hat{K} a_{\frac{N}{2}-1} h\right) \mathrm{FFT}[\psi]]$
		\State $\psi \gets \mathrm{exp}\left(- \ui V b_{\frac{N}{2}-1} h\right) \psi$
		\State \vdots
		\State $\psi \gets \mathrm{FFT}^{-1}[\mathrm{exp}\left(-\ui \hat{K} a_1 h\right) \mathrm{FFT}[\psi]]$
		\State $\psi \gets \mathrm{exp}\left(- \ui V b_1 h \right)\psi$
		\State $t \gets t+h$
		\EndWhile
	\end{algorithmic}
\end{breakablealgorithm}

In the latter, $h$ is the time step, $\mathrm{FFT}$ and $\mathrm{FFT}^{-1}$ denote the Fast 
Fourier Transform and its inverse respectively, $\hat{K}$ is the kinetic 
energy operator in Fourier space, $V$ is the potential, and the values of 
$a_i$ and $b_i$, $i \in \left \{ 1,2,3,4,5,6 \right \}$, are listed in table 
\ref{tSS246}.

\begin{table}[H]
	\centering
	\begin{tabular}{ |c|c|c| } 
		\hline SS2& SS4& SS6\\
		\hline
		$a_1 = \frac{1}{2}$ & $a_1 = \frac{\omega}{2}$ & $a_1 =0.0502627644003922$ \\
		\hline
		$b_1 = 1$ & $b_1 = 1$ & $b_1 =0.148816447901042$ \\
		\hline
		& $a_2 = \frac{1-\omega}{2}$ & $a_2 = 0.413514300428344$ \\
		\hline
		& $b_2 = 1 - 2\omega$ & $b_2 =-0.132385865767784$ \\
		\hline
		& & $a_3 =0.0450798897943977$ \\
		\hline
		& & $b_3 = 0.067307604692185$ \\
		\hline
		& & $a_4 = -0.188054853819569$\\
		\hline
		& & $b_4 = 0.432666402578175$ \\
		\hline
		& & $a_5 = 0.541960678450780$ \\
		\hline
		& & $b_5 = 0.5 - ( b_1 + b_2 + b_3 + b_4)$\\
		\hline
		& & $a_6 = 1 - 2( a_1 + a_2 + a_3 + a_4 + a_5 ) $\\
		\hline
		& & $b_6 = 1 - 2( a_1 + a_2 + a_3 + a_4 + a_5 ) $\\
		\hline
	\end{tabular}
	\caption{Values of the parameters for the Split-Step algorithms. The quantity $\omega$ is given by $\omega=\frac{2+2^{\frac{1}{3}}+2^{-\frac{1}{3}}}{3}$.}
	\label{tSS246}
\end{table}

In the case of SSa, the adaptive splitting algorithm, i.e. the SS4(3), 
both the solutions at the $4^{th}$ and at the $3^{rd}$ order must be 
evaluated. The pseudo-code is described in ALG. (\ref{SSaA}), while the 
coefficients are listed in \ref{tSSa}. In our numerical tests we set 
$\alpha=0.9$, $\beta=3$.

\begin{breakablealgorithm}
	\caption{: SSa}\label{SSaA}
	\begin{algorithmic}[1] 
		\State $t \gets t_0$
		\State $\psi \gets \psi(\boldsymbol{r},t_0)$
		\While{$t < t_{f}$}
		\State $\psi \gets \widetilde{\psi}$
		\State $\psi \gets \mathrm{FFT}^{-1}[\ue^{-\ui \hat{K} a_1 h} \mathrm{FFT}[\psi]]$
		\State $\psi \gets \ue^{- \ui V b_1 h}\psi$
		\State \vdots
		\State $\psi \gets \mathrm{FFT}^{-1}[\ue^{-\ui \hat{K} a_7 h} \mathrm{FFT}[\psi]]$
		\State $\psi \gets \ue^{- \ui V b_7 h}\psi$
		\State $\widetilde{\psi} \gets \mathrm{FFT}^{-1}[\ue^{-\ui \hat{K} \overline{a}_1 h} \mathrm{FFT}[\widetilde{\psi}]]$
		\State $\widetilde{\psi} \gets \ue^{- \ui V \overline{b}_1 h}\widetilde{\psi}$
		\State \vdots
		\State $\widetilde{\psi} \gets \mathrm{FFT}^{-1}[\ue^{-\ui \hat{K} \overline{a}_7 h} \mathrm{FFT}[\widetilde{\psi}]]$
		\State $\widetilde{\psi} \gets \ue^{- \ui V \overline{b}_7 h}\widetilde{\psi}$
		\State $err \gets \sqrt{\frac{ \sum_{i=1}^{N} \left|\psi(\boldsymbol{x}_i,t_n)-\widetilde{\psi}(\boldsymbol{x}_i,t_n)\right|^2 }{\sum_{j=1}^{N}  \left|\psi(\boldsymbol{x}_j,t_n)\right| ^2 }}$
		\If {$err \leq \mathrm{tol}$}
		\State $t \gets t+h$ 
		\Else
		\State $\psi \gets \widetilde{\psi}$
		\EndIf
		\State $h\gets h \min \left \{ \alpha \left(\frac{\mathrm{tol}}{\Delta_n}\right)^{\frac{1}{4}},\beta \right\}$
		\EndWhile
	\end{algorithmic}
\end{breakablealgorithm}

\begin{table}[h]
	\centering
	\begin{tabular}{ |p{0.5cm}|p{4.5cm}|p{0.5cm}|p{4.5cm}|  }
		\hline
		\multicolumn{4}{|c|}{SSa} \\
		\hline
		\multicolumn{2}{|c|}{Order 4} & \multicolumn{2}{|c|}{Order 3}\\
		\hline
		$a_1$& 0& $\Tilde{a}_1$ &0\\
		\hline
		$b_1$& 0.0829844064174052& $\Tilde{b}_1$ &0.0829844064174052\\
		\hline
		$a_2$&0.245298957184271 & $\Tilde{a}_2$ &0.245298957184271\\
		\hline
		$b_2$& 0.3963098014983680& $\Tilde{b}_2$ &0.3963098014983680\\
		\hline
		$a_3$& 0.604872665711080& $\Tilde{a}_3$ &0.604872665711080\\
		\hline
		$b_3$& -0.0390563049223486& $\Tilde{b}_3$ &-0.0390563049223486\\
		\hline
		$a_4$& 0.5 - $(a_2 + a_3)$& $\Tilde{a}_4$ &0.5 - $(a_2 + a_3)$\\
		\hline
		$b_4$& 1. - $2( b_1 + b_2 + b_3 )$& $\Tilde{b}_4$ &1. - $2( b_1 + b_2 + b_3 )$\\
		\hline
		$a_5$& 0.5 - $(a_2 + a_3)$& $\Tilde{a}_5$ &0.3752162693236828\\
		\hline
		$b_5$& -0.0390563049223486& $\Tilde{b}_5$ &0.4463374354420499\\
		\hline
		$a_6$& 0.604872665711080& $\Tilde{a}_6$ &1.4878666594737946\\
		\hline
		$b_6$& 0.3963098014983680& $\Tilde{b}_6$ &-0.0060995324486253\\
		\hline
		$a_7$& 0.245298957184271& $\Tilde{a}_7$ &-1.3630829287974774\\
		\hline
		$b_7$&0.0829844064174052 & $\Tilde{b}_7$ &0\\
		\hline
	\end{tabular}
	\caption{Values of the parameters for the SSa.}
	\label{tSSa}
\end{table}

\end{document}